\documentclass[11pt,twoside,leqno]{amsart}
\usepackage{amssymb,amsbsy,amsmath,amsfonts,amssymb,amscd,epsfig,
times,graphics,color,xypic}

\definecolor{blue}{cmyk}{1.,1.,0.,0.63}
\definecolor{red}{cmyk}{0.,1.,1.,0.63}
\definecolor{green}{cmyk}{1.,0.,1.,0.63}

\usepackage[latin1]{inputenc}
\newcommand{\C}{\mathbb{C}}

\newcommand{\N}{\mathbb{N}}\renewcommand{\P}{\mathbb{P}}

\newtheorem{theorem}{Theorem}
\newtheorem{proposition}{Proposition}
\newtheorem{lemma}{Lemma}

\begin{document}
\parindent 0.75cm

$\:$

\title[
Low pole order frames on jets of the universal hypersurface
]{
Low pole order frames on vertical jets
\\
of the universal hypersurface
}

\author{Jo\"el Merker}

\address{
D\'epartement de Math\'ematiques et Applications, UMR 8553 du CNRS,
\'Ecole Normale Sup\'erieure, 45 rue d'Ulm, F-75230 Paris Cedex 05,
France. \ \
{\tt www.dma.ens.fr/$\sim$merker/}}

\email{merker@dma.ens.fr}

\date{\number\year-\number\month-\number\day}


\begin{abstract}
For low order jets, it is known how to construct
meromorphic frames on the space of the so-called {\sl vertical
$k$-jets}\, $J_{ \sf vert}^k ( \mathcal{ X})$ of the universal
hypersurface $\mathcal{ X} \subset \P^{ n+1} \times \P^{ \frac{ (n+1
+d)!}{ (n+1)!\,\, d!}}$ parametrizing all projective hypersurfaces $X
\subset \P^{ n+1} (\C)$ of degree $d$. In 2004, for $k = n$,
Siu announced that
there exist two constants $c_n \geqslant 1$ and $c_n' \geqslant 1$
such that the twisted tangent bundle:
\[
T_{J_{\rm vert}^n(\mathcal{X})}
\otimes
\mathcal{O}_{\P^{n+1}}
\big(
c_n
\big)
\otimes
\mathcal{O}_{\P^{\scriptscriptstyle{\frac{(n+1+d)!}{(n+1)!\,\,d!}}}}
(c_n')
\]
is generated at every point by its global sections.  In the present
article, we establish this property outside a certain exceptional
algebraic subset $\Sigma \subset J_{\sf vert}^n (\mathcal{ X})$
defined by the vanishing of certain Wronskians, with the {\em
effective} pole order $c_n = \frac{ n^2 + 5n}{ 2}$, thus 
recovering $c_2 = 7$ (Pa\u{u}n), $c_3 = 12$ (Rousseau), 
and with $c_n' = 1$.

Moreover, at the cost of raising $c_n$ up to $c_n = n^2 + 2n$, the
same generation property holds outside the smaller set $\widetilde{
\Sigma} \subset \Sigma \subset J_{ \sf vert}^n ( \mathcal{ X})$ which
is defined by the vanishing of all first order jets. Applications to
{\em weak} (with $\Sigma$) and to {\em strong} (with $\widetilde{
\Sigma}$) algebraic degeneracy of entire holomorphic curves 
$\C \to X$ are 
upcoming.
\end{abstract}

\maketitle

\vspace{-0.75cm}

\begin{center}
\begin{minipage}[t]{10cm}
\baselineskip =0.35cm
{\scriptsize

\centerline{\bf Table of contents}

\smallskip

{\bf 1.~Introduction \dotfill \pageref{Section-1}.}

{\bf 2.~Universal hypersurface and vertical jets \dotfill
\pageref{Section-2}.}

{\bf 3.~First package of coefficient vector fields 
\dotfill \pageref{Section-3}.}

{\bf 4.~Second package of jet, coordinate vector fields 
\dotfill \pageref{Section-4}.}

{\bf 5.~Appendix: a determinantal identity
\dotfill \pageref{Section-5}.}

\smallskip

}\end{minipage}
\end{center}

\section*{ \S1.~Introduction}
\label{Section-1}

The Kobayashi hyperbolicity conjecture (1970), in optimal degree and
taking account of Brody's theorem (1978), expects that all entire
holomorphic curves $f : \C \to X$ into a complex projective
(algebraic, smooth) hypersurface $X \subset \P^{ n+1}$ must be
constant if $\deg X \geqslant 2n+1$, provided $X$ is generic.

In 2004, Siu~\cite{ siu2004} announced a strategy of proof, valid in
arbitrary dimensions for (extremely) high 
(noneffective) degrees $d >\!\!> n$. Two
major techniques are used.

Inspired by Bloch's ideas, one looks firstly for global sections of
the Green-Griffiths bundle $E_{ k,m}^{ GG} T_X^*$ of jet differentials
of order $k$ and weighted degree $m$ ({\em cf.} \cite{ gg1980}), which
vanish on some ample divisor; an Ahlfors-Schwarz-type theorem then
forces every entire curve $f: \C \to X$ to satisfy the corresponding
differential equation (\cite{ dem1997}), a first step toward algebraic
degeneracy.  In 1997, Demailly introduced a refined subbundle $E_{
k,m} T_X^*$ having better positivity properties which consists of jet
differentials that are invariant under (local) reparametrizations of
the source $\C$. In dimension $n = 3$ for jets of order $k = 3$,
Rousseau (\cite{ rou2006a}) completely described the algebraic
structure of $E_{ k,m} T_X^*$ in its fibers, decomposed it in direct
sums of Schur bundles $\Gamma^{ (\lambda_1, \lambda_2, \lambda_3)}
T_X^*$, computed its Euler characteristic $\chi \big( X, \, E_{ k,m}
T_X^* \big)$, majorated from above $h^2 \big( X, \, E_{ k,m} T_X^*
\big)$ ({\em see} \cite{ rou2006b}), and established existence of
global algebraic differential equations in degree $d \geqslant 97$.

In~\cite{ mer2008a, mer2008b}, one finds a {\em complete algorithm} to
generate all Demailly-Semple 
invariants in arbitrary dimension $n \geqslant
2$. In particular, for $n = k = 4$, there are 16 fundamental, mutually
independent bi-invariant polynomials generating the Demailly-Semple
(unipotent-invariant) algebra sharing 41 (groebnerized) syzygies, and
one deduces by polarization
that the algebra of all invariants for $n = k =
4$ is generated by 2835 polynomials. Nonconstant entire holomorphic
curves valued in an algebraic 3-fold (resp. 4-fold) $X^3 \subset \P^4
(\C)$ (resp. $X^4 \subset \P^5 (\C)$) of degree $d$ satisfy (\cite{
mer2008b}) global differential equations as soon as $d \geqslant 72$
(resp. $d \geqslant 259$).

In~\cite{ div2007}, for dimensions $n = 2, 3, 4, 5$ and for jet orders
$k = 3, 4, 5, \underline{ 5}$, resp., it is shown that asymptotically
as $m \to \infty$:
\[
H^0\big(X,\,E_{k,m}T_X^*\otimes A^{-1}\big)
\neq
0,
\]
in degrees $d\geqslant 16$, $74$, $298$, $1222$ resp., where $A \to X$
is any auxiliary ample line bundle.  But it is also shown that $H^0
\big( X, \, E_{ k,m} T_X^* \big) = 0$, for all jet orders $k \leqslant
\dim X - 1$, generalizing a theorem of Rousseau (\cite{ rou2006b}) in
dimension $3$. Furthermore, 
for jet order $k$ equal to the dimension $n$,
with $n$ arbitrary, Diverio shows in~\cite{ div2008} that there exists
an integer $\delta_n >\!\!> n$ (up to now not effective) insuring
existence of global sections of $E_{ n,m} T_X^* \otimes A^{ -1}$ in
degree $d \geqslant \delta_n$.

The second technique, initiated by Clemens~\cite{ cle1986}, Ein~\cite{
ein1988}, Voisin~\cite{ voi1996} and pushed further by Siu~\cite{
siu2004}, Pa\u{u}n~\cite{ pau2008}, Rousseau~\cite{ rou2007a},
consists in constructing meromorphic frames on the space of the
so-called {\sl vertical $k$-jets} $J_{\sf vert}^k ( \mathcal{ X})$ in
the universal hypersurface $\mathcal{ X} \subset \P^{ n+1} \times \P^{
{\scriptscriptstyle{ \frac{ (n+1 +d)!}{ (n+1)!\,\, d!}}}}$
parametrizing all $X \subset \P^{ n+1}$ of degree $d$, so as to
produce, by frame differentiations, enough {\em independent}\,
algebraic differential equations from just one global section of $E_{
k,m} T_X^* \otimes A^{ -1}$.

In~\cite{ siu2004}, p.~557, Siu
announced that, for $k = n$, there exist two constants $c_n \geqslant
1$ and $c_n' \geqslant 1$ such that the twisted tangent bundle:
\[
T_{J_{\rm vert}^n(\mathcal{X})}
\otimes
\mathcal{O}_{\P^{n+1}}
\big(
c_n
\big)
\otimes
\mathcal{O}_{\P^{\scriptscriptstyle{\frac{(n+1+d)!}{(n+1)!\,\,d!}}}}
(c_n')
\]
is generated at every point by its global sections (frame
property). In the present article, we establish this property outside
a certain exceptional algebraic subset $\Sigma \subset J_{\sf vert}^n
(\mathcal{ X})$ defined by the vanishing of certain Wronskians, with
the effecive pole order $c_n = \frac{ n^2 + 5n}{ 2}$, recovering $c_2
= 7$ (Pa\u{u}n \cite{ pau2008}), $c_3 = 12$ (Rousseau \cite{
rou2007a}), and with $c_n' = 1$.

Moreover, at the cost of raising $c_n$ up to $c_n = n^2 + 2n$, the same
generation property holds outside the smaller set $\widetilde{ \Sigma}
\subset \Sigma$ defined by the vanishing of all first order
jets. Applications to {\em weak} (with $\Sigma$) and to {\em strong}
(with $\widetilde{ \Sigma}$) algebraic degeneracy of entire
holomorphic curves are given in \cite{ dmr2008}, following Rousseau's
Schur bundle decomposition strategy in 
dimension $n = 4$, and also in higher dimensions,
thanks to Diverio's use (\cite{ div2007, div2008}) of the algebraic
version of Demailly's Morse inequalities due to Trapani (\cite{
tra1995}).

\subsection*{ Acknowledgments} 
During the author's stay at the Mittag-Leffler Institute (1--21 April
2008), Yum-Tong Siu and Mihai Pa\u{u}n have provided helpful oral
explanations of~\cite{ siu2004, pau2008}.  The stronger property of
generation at every point of $J_{\sf vert}^n ( \mathcal{ X})
\backslash \widetilde{ \Sigma}$ was obtained thanks to fruitful
exchanges joint with Simone Diverio and Erwan Rousseau during the {\em
Workshop Complex Hyperbolic Geometry and Related Topics} at the
Fields Institute, Toronto, Canada, 17-21 November 2008.

\section*{ \S2.~Universal hypersurface and vertical jets}
\label{Section-2}

\subsection*{ Representation in coordinates}
Consider the {\sl universal hypersurface}\, $\mathcal{ X} \subset \P^{
n + 1} \times \P^{ \frac{ (n+1+d)!}{(n+1)!  \,\,d!}}$ parametrizing
all complex $n$-dimensional algebraic hypersurfaces of fixed degree $d
\geqslant 1$ in $\P^{ n + 1}$ which is defined, in two collections of
homogeneous coordinates:
\[
\aligned
\left[Z\right]
&
=
[Z_0:Z_1:\cdots:Z_n:Z_{n+1}]
\in\P^{n+1}
\\
[A]
&
=
\big[(A_\alpha)_{\alpha\in\N^{n+2},\,\vert\alpha\vert=d}\big]
\in\P^{\frac{(n+1+d)!}{(n+1)!\,\,d!}},
\endaligned
\]
as the zero-set locus:
\[
\mathcal{X}\ :
\ \ \ \ \ \ \ \ \
0
=
\sum_{\alpha\in\N^{n+2}\atop 
\vert\alpha\vert=d}\,
A_\alpha\,Z^\alpha
\]
of the general homogeneous degree $d$ polynomial. Here of course, a
multiindex $\alpha = ( \alpha_0, \alpha_1, \dots, \alpha_{ n+1}) \in
\N^{ n+2}$ has {\sl length} defined by $\vert \alpha \vert := \alpha_0 +
\alpha_1 + \cdots + \alpha_{ n+1}$ and we abbreviate $Z^\alpha =
Z_0^{\alpha_0} Z_1^{\alpha_1} \cdots Z_{n+1}^{ \alpha_{n+1}}$.

Our goal is to perform, for jets of order $\kappa$ equal to the
dimension $n$ of hypersurfaces $\mathcal{ X} ( A) \subset \P^{ n+1}$,
a construction of meromorphic vector fields on the space of jets of
holomorphic discs (or entire maps) valued in $\mathcal{ X}$ which was
initiated by Clemens~\cite{ cle1986}, Ein~\cite{ ein1988},
Voisin~\cite{ voi1996} for $\kappa = 1$, $n \geqslant 1$, then
announced for higher $\kappa$'s by Siu~\cite{ siu2004} and recently
detailed by Pa\u{u}n~\cite{ pau2008} for $n = \kappa = 2$ and by
Rousseau~\cite{ rou2007a} for $n = \kappa = 3$.
For general $\kappa = n$, a concise book-keeping of indices appears to
be available here. 

As in~\cite{ pau2008, rou2007a}, we shall mainly work in {\em
in}homogeneous coordinates on the chart $\{ Z_0 \neq 0 \} \times \{
A_{ 0d0\cdots 0} \neq 0\}$, a copy of $\C^{ n + 1} \times \C^{ \frac{
( n +1 + d)!}{ (n+1)! \, \, d!}}$.  Dividing by $(Z_0)^d$ and by $A_{
0 d 0 \cdots 0}$, and setting $z_i := Z_i / Z_0$, the equation of
$\mathcal{ X}$ then transfers to:
\[
\mathcal{X}_0\ :
\ \ \ \ \ \ \ \
0
=
z_1^d
+
\sum_{\alpha\in\N^{n+1}\atop\vert\alpha\vert\leqslant d,\,
\alpha_1<d}\,\,
a_\alpha\,z^\alpha,
\]
with new coefficients $a_{ \alpha_1 \cdots \alpha_{ n+1}} := \frac{
A_{ \alpha_0 \alpha_1 \cdots \alpha_{ n+1}}}{ A_{ 0 d 0 \cdots 0}}$ in
which $\alpha_0 := d - \alpha_1 - \cdots - \alpha_{ n+1}$. By
convention, we shall set $a_{ d 0 \cdots 0} = 1$.

In view of applications to the Green-Griffiths algebraic degeneracy
conjecture ($d \geqslant n + 3$) or to the Kobayashi hyperbolicity
conjecture ($d \geqslant 2 n + 1$), it will, without
loss of generality, be assumed that $d > n$ throughout.

\subsection*{ Defining equations for the space of vertical jets}
To settle Kobayashi hyperbolicity or Green-Griffiths algebraic
degeneracy, the strategy initiated by Bloch and pursued by
Green-Griffiths~\cite{ gg1980}, Siu~\cite{ siu2004}, Demailly~\cite{
dem1997} consists in producing enough (global, algebraic) differential
equations that every entire map $\C \ni \zeta \longmapsto \big( z_1 (
\zeta), \dots, z_{ n+1} ( \zeta) \big)$ valued in an algebraic variety
$\mathcal{ X} (A)$ for (very) generic fixed coefficients $A_\alpha$
should satisfy. Accordingly, if one introduces independent coordinates 
corresponding to derivatives with respect to $\zeta$:
\[
\Big(z_i,\,a_\alpha,\,z_{j_1}',\,z_{j_2}'',\dots,
z_{j_\kappa}^{(\kappa)}\Big)
\in\C^{n+1}\times\C^{\frac{(n+1+d)!}{(n+1)!\,\,d!}}\times
\underbrace{\C^{n+1}\times
\C^{n+1}\times\cdots\times\C^{n+1}}_{\kappa\,\,\text{\rm times}}\,,
\]
the manifold of $\kappa$-jets of such entire maps has equations
obtained by just formally differentiating the monomials $z^\alpha$
with respect to the variable $\zeta \in \C$, the $a_\alpha$ being
constant. The basic chain rule yields the first five equations, up to
$\kappa = 4$:
\[
\footnotesize
\aligned
0
&
=
\sum_{\alpha\in\N^{n+1}\atop
\vert\alpha\vert\leqslant d,\,a_{d0\cdots 0}=1}\,
a_\alpha\,z^\alpha
\\
0
&
=
\sum_\alpha\,a_\alpha
\bigg(
\sum_{j_1}\,\frac{\partial (z^\alpha)}{\partial z_{j_1}}\,z_{j_1}'
\bigg)
\\
0
&
=
\sum_{\alpha}\,a_\alpha
\bigg(
\sum_{j_1}\,\frac{\partial (z^\alpha)}{\partial z_{j_1}}\,z_{j_1}''
+
\sum_{j_1,\,j_2}\,
\frac{\partial^2 (z^\alpha)}{\partial z_{j_1}\partial z_{j_2}}\,
z_{j_1}'z_{j_2}'
\bigg)
\\
0
&
=
\sum_\alpha\,a_\alpha\,
\bigg(
\sum_{j_1}\,\frac{\partial (z^\alpha)}{\partial z_{j_1}}\,z_{j_1}'''
+
\sum_{j_1,\,j_2}\,
\frac{\partial^2(z^\alpha)}{\partial z_{j_1}\partial z_{j_2}}\,
3\,z_{j_1}'z_{j_2}''
+
\sum_{j_1,\,j_2,\,j_3}\,
\frac{\partial^3(z^\alpha)}{\partial z_{j_1}\partial z_{j_2}\partial
z_{j_3}}\,
z_{j_1}'z_{j_2}'z_{j_3}'
\bigg)
\\
0
&
=
\sum_\alpha\,a_\alpha
\bigg(
\sum_{j_1}\,
\frac{\partial (z^\alpha)}{\partial z_{j_1}}\,z_{j_1}''''
+
\sum_{j_1,\,j_2}\,
\frac{\partial^2(z^\alpha)}{\partial z_{j_1}\partial z_{j_2}}
\big(
4\,z_{j_1}'z_{j_2}'''
+
3\,z_{j_1}''z_{j_2}''
\big)
+
\\
&
\ \ \ \ \ \
+
\sum_{j_1,\,j_2,\,j_3}\,
\frac{\partial^3(z^\alpha)}{\partial z_{j_1}
\partial z_{j_2}\partial z_{j_3}}\,
6\,z_{j_1}'z_{j_2}'z_{j_3}''
+
\sum_{j_1,\,j_2,\,j_3,\,j_4}\,
\frac{\partial^4(z^\alpha)}{\partial z_{j_1}\partial
z_{j_2}\partial z_{j_3}\partial z_{j_4}}\,
z_{j_1}'z_{j_2}'z_{j_3}'z_{j_4}'
\bigg),
\endaligned
\]
on understanding that $a_{ d 0 \cdots 0} = 1$ and that all summations
$\sum_{ j_1}$, $\sum_{ j_1, \, j_2}$ {\em etc.}  are performed for the
indices $j_i$ running from $1$ to $n+1$. Equivalently, this
submanifold of $\C^{ n+1} \times \C^{ \frac{ ( n + 1+ d)!}{ (n+1)!  \,
\, d!}}  \times \C^{\kappa (n + 1)}$ may be be defined as the
submanifold of the full $\kappa$-jet manifold $J^\kappa ( \C, \,
\mathcal{ X}_0)$ consisting of only the jets tangent to the fibers of
the projection $\mathcal{ X}_0 \to \P^{ \frac{ ( n + 1+ d)! }{ (n+1)!
\, \, d!}}$ onto the second factor.  They are called {\sl vertical
jets} in~\cite{ siu2004, pau2008, rou2007a} and will be denoted by
$J_{\sf vert}^\kappa ( \mathcal{ X}_0)$.

Formally differentiating any polynomial in the jet variables amounts
to applying the {\sl total differentiation operator}\:\!:
\[
{\sf D}
(\bullet)
:=
\sum_{\lambda\in\N}\,\,
\sum_{k=1}^{n+1}\,
\frac{\partial (\bullet)}{\partial z_k^{(\lambda)}}\,
\cdot\,
z_k^{(\lambda+1)},
\]
and above, it is clear that each next equation is obtained from the
previous one by applying ${\sf D}$ to it so that, for jets of
arbitrary order $\kappa$ up to $\kappa$ equal to the dimension $n$,
the $( n+1 )$ defining equations of $J_{\sf vert}^n ( \mathcal{ X}_0)$
happen to be:
\[
0
=
\sum_\alpha\,a_\alpha z^\alpha
=
{\sf D}
\Big(
\sum_\alpha\,a_\alpha z^\alpha
\Big)
=\cdots=
{\sf D}^n
\Big(
\sum_\alpha\,a_\alpha z^\alpha
\Big).
\]
Then a suitable multivariate version of the classical Fa� di Bruno
formula provides a {\sl closed, explicit formula} for all such
equations.

\def\thelemma{1}\begin{lemma}
{\rm (\cite{cs1996, mer2005})} The $(n + 1)$ defining equations of
$J_{\sf vert}^n ( \mathcal{ X}_0)$ write as follows, where $\kappa =
0, 1, 2, \dots, n${\rm :}
\[
\footnotesize
\aligned
0
&
=\!\!\!\!\!\!
\sum_{\alpha\in\N^{n+1}\atop
\vert\alpha\vert\leqslant d,\,a_{d0\cdots0}=1}\!\!\!\!\!
a_{\alpha}\,\,\,
\sum_{e=1}^\kappa\,
\sum_{1\leqslant\lambda_1<\cdots<\lambda_e\leqslant\kappa}\,
\sum_{\mu_1\geqslant 1,\dots,\mu_e\geqslant 1}\,
\sum_{\mu_1\lambda_1+\cdots+\mu_e\lambda_e=\kappa}\,
\frac{\kappa!}{(\lambda_1!)^{\mu_1}\mu_1!\cdots
(\lambda_e!)^{\mu_e}\mu_e!}\,
\\
&
\sum_{j_1^1,\dots,j_{\mu_1}^1=1}^{n+1}\cdots
\sum_{j_1^e,\dots,j_{\mu_e}^e=1}^{n+1}
\frac{\partial^{\mu_1+\cdots+\mu_e}\big(z^\alpha\big)}{\partial
z_{j_1^1}\cdots\partial z_{j_{\mu_1}^1}\cdots
\partial z_{j_1^e}\cdots\partial z_{j_{\mu_e}^e}}\,\,
z_{j_1^1}^{(\lambda_1)}\cdots z_{j_{\mu_1}^1}^{(\lambda_1)}
\cdots
z_{j_1^e}^{(\lambda_e)}\cdots z_{j_{\mu_e}^e}^{(\lambda_e)}.
\endaligned
\]
\end{lemma}

To read this general formula with the help of the formulas specialized
above, we comment it backwards from its end.

The general monomial $\prod \, z_{\bullet}^{( \lambda_1)} \, \prod
z_{\bullet}^{ ( \lambda_2)} \cdots \prod\, z_{\bullet }^{ (
\lambda_e)}$ in the jet variables gathers derivatives of increasing
orders $\lambda_1 < \lambda_2 < \cdots < \lambda_e$, with $\mu_1,
\mu_2, \dots, \mu_e$ counting their respective numbers.  Then each
monomial $z^\alpha$ is subjected to a partial derivative of order
$\mu_1 + \mu_2 + \cdots + \mu_e$, the total number of $z_j^{ (
\lambda_i)}$ in the monomial in question. Since there are $n + 1$
variables $z_i$, the dots in the $z_{ \bullet}^{ ( \lambda_i)}$ should
receive indices, and in fact, there appear general sums $\sum_{ j_1^i,
\dots, j_{ \mu_i}^i = 1}^{ n +1}$ over {\em all possible} such
indices. Notice that these observations are confirmed by the
formulas developed above up to $\kappa = 4$.

In the sequel, we will in fact not need all the information of such a
precise, explicit formula, but it will suffice to know that, among the
$(n+1)$ defining equations, the equation numbered $\kappa$ is a
certain finite sum with certain integer coefficients of terms of 
the form:
\[
\boxed{
\sum_{\beta\in\N^{n+1}\atop
\vert\beta\vert\leqslant d,\,a_{d0\cdots0}=1}\,a_{\beta}\,
\bigg(
\sum_{j_1,\dots,j_e=1}^{n+1}\,\,
\frac{\partial^e(z^\beta)}{\partial z_{j_1}\cdots
\partial z_{j_e}}
\cdot
z_{j_1}^{(\nu_1)}\cdots
z_{j_e}^{(\nu_e)}
\bigg)}\,,
\]
where the derivative orders $\nu_i \geqslant 1$ of the jet monomial
$z_{ j_1 }^{ ( \nu_1 ) } \cdots z_{ j_e }^{ ( \nu_e ) }$ are
nondecreasing and where $\nu_1 + \cdots + \nu_e = \kappa$.
The reader unacquainted with the Fa� di Bruno combinatorics
could readily prove this less informative representation by 
reasoning inductively on $\kappa$.

\subsection*{ Frames and generation by global sections}
Now, a globally defined vector field on the ambient space $\C^{ n + 1}
\times \C^{ \frac{ (n+1+d)!}{ (n+1)! \,\, d!}} \times \C^{ n ( n+1)}$
writes under the general form:
\[
\footnotesize
\aligned
{\sf T}
=
\sum_{i=1}^{n+1}\,
{\sf Z}_i\,\frac{\partial}{\partial z_i}
+
\sum_{\alpha\in\N^{n+1}\atop\vert\alpha\vert\leqslant d,\,
\alpha_1<d}\,
{\sf A}_\alpha\,
\frac{\partial}{\partial a_\alpha}
+
\sum_{k=1}^{n+1}\,{\sf Z}_k'\,\frac{\partial}{\partial z_k'}
+
\sum_{k=1}^{n+1}\,{\sf Z}_k''\,\frac{\partial}{\partial z_k''}
+\cdots+
\sum_{k=1}^{n+1}\,{\sf Z}_k^{(n)}\,\frac{\partial}{\partial z_k^{(n)}}.
\endaligned
\]
We shall seek vector fields of this form which should extend
meromorphically to the full space of vertical jets and which should
make a spanning frame of vectors tangent to $J_{ \sf vert}^n (
\mathcal{ X})$ at almost every point, say outside a certain ``bad''
set. After twisting by $(\bullet) \otimes \mathcal{ O}_{ \P^{ n+1}} (
c_n ) \otimes \mathcal{O}_{ \P^{ \scriptscriptstyle{ \frac{ (n+1+d)!}{
(n+1)! \,\,d! }}}} (c_n')$ for some two suitable constants $c_n
\geqslant 1$ and $c_n' \geqslant 1$, one may in fact erase the
appearing poles of the
meromorphic coefficients, 
so that one may speak of global {\em holomorphic} sections
instead of meromorphic sections.

\def\thetheorem{\!}\begin{theorem}
Let $\widetilde{
\Sigma}$ be the closure, in $J_{\sf vert}^n ( \mathcal{ X})$, of
the Zariski closed subset of the space $J_{\sf vert}^n ( \mathcal{
X}_0)$ of vertical {\rm affine} jets defined by requiring that all
{\em first order} jet vanish:
\[
\widetilde{\Sigma}_0
:=
\Big\{
\big(z_i,a_\alpha,z_{j_1}',\dots,z_{j_n}^{(n)}\big):\,\,
z_1'=z_2'=\cdots=z_{n+1}'=0
\Big\},
\]
so that in any other standard affine chart $(t_0, \dots, t_{ \upsilon
  -1}, t_{ \upsilon + 1}, \dots, t_{ n+1}) \in \C^{ n+1}$ on $\P^{
  n+1} ( \C)$, the representation of $\widetilde{ \Sigma}$ is yielded
by exactly the same equations $0 = t_0' = \cdots = t_{ \upsilon - 1} '
= t_{ \upsilon + 1}' = \cdots = t_{ n+1}'$. Then the following two
properties hold true.

\begin{itemize}

\smallskip\item[$\bullet$]
$J^n_{ \sf vert} ( \mathcal{ X}) \backslash \Sigma$ is smooth of pure
codimension equal to $n+1$ at every point, namely, it is of 
pure dimension equal to:
\[
\aligned
j_n^d
:=
&\
n+1+
{\textstyle{\frac{(n+1+d)!}{(n+1)!\,d!}}}
+n(n+1)
-(n+1)
\\
=
&\
{\textstyle{\frac{(n+1+d)!}{(n+1)!\,d!}}}
+
n(n+1).
\endaligned
\]

\smallskip\item[$\bullet$]
The twisted tangent bundle:
\[
T_{J_{\rm vert}^n(\mathcal{X})}
\otimes
\mathcal{O}_{\P^{n+1}}
\big(
n^2+2n
\big)
\otimes
\mathcal{O}_{\P^{\scriptscriptstyle{\frac{(n+1+d)!}{(n+1)!\,\,d!}}}}
(1)
\]
is generated by its global sections on $J_{\sf vert}^n ( \mathcal{ X})
\big \backslash \widetilde{ \Sigma}$, that is to say: at every point
$p^{ [ n]} \in J_{\sf vert}^n ( \mathcal{ X}) \big \backslash
\widetilde{ \Sigma}$ not lying in $\widetilde{ \Sigma}$, one may find
$j_n^d$ {\em global} sections ${\sf T}_1, \dots, {\sf T}_{ j_n^d}$
over $X$ of this twisted tangent bundle such that:
\[
\C{\sf T}_1(p)\oplus\cdots\oplus\C{\sf T}_{j_n^d}(p)
=
T_{J_{\sf vert}^n(\mathcal{X}),\,p}.
\]
\end{itemize}
\end{theorem}

\subsection*{ Comments about applications}
The simplicity of the defining equations of the avoided ``bad''
set $\widetilde{
\Sigma} = \{ z_i' = 0 \}$ has considerable advantages in the study of
Green-Griffiths algebraic degeneracy of entire holomorphic curves $f :
\C \to X$.

Indeed, by employing jet differentials, one shows in a first moment
that the $n$-jet $j^nf$ of any such an $f$ must satisfy\footnote{\,
{\em see} \cite{ gg1980, dem1997, siu2004, rou2006a, rou2006b, div2007}; we only summarize very briefly the ideas here.}  at least
one nontrivial global algebraic differential equation $P ( j^n f ) =
0$.  Then in a second moment, following Siu's strategy ({\em see}
\cite{ siu2004, pau2008, rou2007a, dmr2008}) which consists in
applying some well chosen multi-derivations $({\sf T}_1)^{ \nu_1}
\cdots ({\sf T}_{ j_n^d})^{ \nu_{ j_n^d}}$ to $P ( j^n f) = 0$ so as
to get sufficiently many {\em supplementary} differential equations,
one comes down to distinguishing two cases:

\begin{itemize}

\smallskip\item[$\square$\,\,]
either $j^n f ( \C) \not \subset 
\widetilde{ \Sigma}$; in this first case, one is
then able to show (\cite{ rou2007a, dmr2008}) that $f ( \C)$ is
contained in a certain proper algebraic subvariety $Y \subsetneqq X$
which is independent of $f$, and this yields {\em strong} algebraic
degeneracy;

\smallskip\item[$\square$\,\,]
or else $j^n f ( \C) \subset \widetilde{ \Sigma}$ fully; in this
second case, one cannot apply any derivation ${\sf T}_1, \dots,
{\sf T}_{ j_n^d}$, but then the condition $j^n f ( \C) \subset \Sigma$
simply reads $0 \equiv f_1' (\zeta) \equiv f_2 ' (\zeta) \equiv \cdots
\equiv f_n' (\zeta)$, hence $f$ is {\em constant} and strong
degeneracy again holds {\em gratuitously}.

\end{itemize}\smallskip

Quite differently, in~\cite{ pau2008, rou2007a} and in a preliminary
version of the present article as well, the ``bad'' set $\Sigma$ that
one had to avoid was substantially larger than $\widetilde{
\Sigma}$. Then as a consequence in these references, the condition
$j^n f \subset \Sigma$ in the second case above only meant that $j^n
f$ was contained in the intersection of $X$ with some
one-codimensional linear subspace $H$ of $\P^{ n+1} ( \C)$ which in
general depended upon $f$, so that only {\em weak} algebraic
degeneracy of $f ( \C)$ could be deduced\footnote{\, A more careful
inspection shows that in fact, $H$ is two-codimensional
(Simone Diverio). }. Here is
the weaker statement which we generalize in arbitrary dimension $n
\geqslant 2$.

\def\thetheorem{\!\!'}\begin{theorem}
Let  $\Sigma$ be the closure, in $J_{\sf vert}^n ( \mathcal{ X})$, of
the Zariski closed subset of the space $J_{\sf vert}^n (
\mathcal{ X}_0)$ of vertical {\rm affine} jets defined by 
requiring that all $n \times n$ Wronskians vanish:
\[
\aligned
\Sigma_0
:=
\Big\{
\big(z_i,a_\alpha,z_{j_1}',\dots,z_{j_n}^{(n)}\big):\,\,
0
&
=
\det\,\big(z_i^{(\lambda_j)}\big)_{1\leqslant i\leqslant n+1}^{
1\leqslant j\leqslant n}
\\
&
\text{for all}\,\,\lambda_1,\dots,\lambda_n\,\,
\text{with}\,\,
1\leqslant\lambda_j\leqslant n
\Big\}.
\endaligned
\]  
Then the twisted tangent bundle:
\[
T_{J_{\rm vert}^n(\mathcal{X})}
\otimes
\mathcal{O}_{\P^{n+1}}
\big(
{\textstyle{\frac{n^2+5n}{2}}}
\big)
\otimes
\mathcal{O}_{\P^{\scriptscriptstyle{\frac{(n+1+d)!}{(n+1)!\,\,d!}}}}
(1)
\]
is generated by its global sections at every 
point of $J_{\sf vert}^n ( \mathcal{ X}) \big \backslash \Sigma$.
\end{theorem}

Notice that the twisting order $\frac{ n^2 + 5n }{ 2}$ along the
$z$-direction is smaller than the one $n^2 + 2n$ of the preceding
theorem: a certain price has to be ``paid'' in order to
shrink the ``bad'' set, and to thereby gain strong degeneracy.

As said, one may verify that the vanishing of all $n \times n$ minors
of the $n \times (n+1)$ Wronskian-like matrix $\big( f_j^{ (\lambda)}
(\zeta) \big)_{ 1 \leqslant j \leqslant n+1}^{ 1 \leqslant \lambda
\leqslant n}$ implies that the components $f_1 (\zeta), \dots, f_{
n+1} ( \zeta)$ satisfy at least two linearly independent linear
relations:
\[
0
\equiv
{\textstyle{\sum_{i=1}^{n+1}}}\,a_i\,f_i(\zeta)
\equiv
{\textstyle{\sum_{i=1}^{n+1}}}\,b_i\,f_i(\zeta), 
\]
for
$\zeta \in \C$, with no universal control on the coefficients 
$a_i, b_i$. 

\smallskip

Before proceeding to establishing the two theorems, let us check that
the set $\Sigma$ is represented by the same kind of equations $0 =
t_1' = \cdots = t_{ \upsilon - 1}' = t_{ \upsilon + 1}' = \cdots = t_{
n+1}'$ in any other standard chart $\{ Z_\upsilon \neq 0\}$ on $\P^{
n+1} ( \C)$ in which the affine coordinates are defined just by:
\[
t_0
=
{\textstyle{\frac{Z_0}{Z_\upsilon}}},
\dots,
t_{\upsilon-1}
=
{\textstyle{\frac{Z_{\upsilon-1}}{Z_\upsilon}}},\,\,
t_{\upsilon+1}
=
{\textstyle{\frac{Z_{\upsilon+1}}{Z_\upsilon}}},
\dots,
t_{n+1}
=
{\textstyle{\frac{Z_{n+1}}{Z_\upsilon}}}.
\]
Indeed, coming back to the definition $z_i = \frac{ Z_i}{ Z_0}$ of the
$z_i$, $i=1, \dots, n+1$, the change of chart $\{ Z_0 \neq 0 \} \to \{
Z_\upsilon \neq 0 \}$ is given by the well known basic formulas:
\[
t_0
=
{\textstyle{\frac{1}{z_\upsilon}}},
\dots,
t_{\upsilon-1}
=
{\textstyle{\frac{z_{\upsilon-1}}{z_\upsilon}}},\,\,
t_{\upsilon+1}
=
{\textstyle{\frac{z_{\upsilon+1}}{z_\upsilon}}},
\dots,
t_{n+1}
=
{\textstyle{\frac{z_{n+1}}{z_\upsilon}}},
\]
whence by differentiating the right-hand sides as if
they virtually depended upon a variable $\zeta \in \C$, 
we get the transformation rules for the first order jets:
\[
t_0'
=
-{\textstyle{\frac{z_\upsilon'}{z_\upsilon^2}}},
\dots,
t_{\upsilon-1}'
=
{\textstyle{\frac{z_{\upsilon-1}'}{z_\upsilon}}}
-
{\textstyle{\frac{z_{\upsilon-1}z_\upsilon'}{z_\upsilon^2}}},\,\,
t_{\upsilon+1}'
=
{\textstyle{\frac{z_{\upsilon+1}'}{z_\upsilon}}}
-
{\textstyle{\frac{z_{\upsilon+1}z_\upsilon'}{z_\upsilon^2}}},
\dots,
t_{n+1}'
=
{\textstyle{\frac{z_{n+1}'}{z_\upsilon}}}
-
{\textstyle{\frac{z_{n+1}z_\upsilon'}{z_\upsilon^2}}}.
\]
Then visibly, the two representations $\{ 0 = z_1' = \cdots = z_{
n+1}' \}$ and $\{ 0 = t_0' = \cdots = t_{ \upsilon -1} ' = t_{
\upsilon + 1} ' = \cdots = t_{ n+1} ' \}$ of the set $\Sigma$ coincide
coherently on the intersection $\{ Z_0 \neq 0 \} \cap \{ Z_\upsilon
\neq 0 \}$ of the two affine charts. One may verify that $\widetilde{
\Sigma}$ also enjoys a similar invariance property.

\smallskip

\subsection*{ Organization}
The remainder of the paper is entirely devoted to the proof of the
first theorem. When necessary, we shall briefly indicate which
mild modifications suffice to gain the second theorem at
the same time.

\section*{\S3.~First package of coefficient vector fields}
\label{Section-3}

\subsection*{ First family of global sections}
We begin by seeking tangent vector fields globally defined over $\C^{
n+1} \times \C^{ \frac{ (n+1+d)!}{ (n+1)! \, \, d!}} \times \C^{
n(n+1)}$ of the specific, short form:
\[
{\sf T}
=
\sum_{\vert\alpha\vert\leqslant n}\,
{\sf A}_\alpha\,
\frac{\partial}{\partial a_\alpha}
\]
in the space of only the coefficient variables $a_\alpha$, up to
length $n$. Afterwards, we shall deal with $\sum_{ n \leqslant \vert
\alpha \vert \leqslant d \atop \alpha_1 < d}\, {\sf A}_\alpha \,
\frac{ \partial }{\partial a_\alpha}$, and in Section~4, the remaining
directions $\partial \big/ \partial z_i$ and $\partial \big/ \partial
z_j^{ (\lambda)}$ will complete the sought generating tangent
vector fields.

Any arbitrary point $p^{ [n]} \in J_{\sf vert}^n ( \mathcal{ X}_0)$
not in $\widetilde{ \Sigma}$ lies in at least one of
the open sets $\{ z_i'
\neq 0\}$. Fixing such an index $i$ with
$1 \leqslant i \leqslant n+1$, we shall construct a
collection of vector fields of the above form that are defined in $\{
z_i' \neq 0\}$ and that extend meromorphically to $J_{\sf vert}^n (
\mathcal{ X})$.  To this aim, let us rewrite the defining equations of
$J_{\sf vert}^n ( \mathcal{ X}_0)$ under the following convenient
form, in which we denote by $\epsilon_i = (0, \dots, 1, \dots, 0)$ the
$i$-th basic multiindex having $1$ at the $i$-th place and $0$
elsewhere, whence $n \epsilon_i = (0, \dots, n, \dots, 0)$:
\[
\left\{
\aligned
0
&
=
a_0
+
a_{\epsilon_i}z_i
+\cdots\cdots\cdots\cdot\cdot+
a_{n\epsilon_i}z_i^n
+
\sum_{\beta\neq\epsilon_i,\dots,n\epsilon_i\atop
1\leqslant\vert\beta\vert\leqslant d}\,
a_\beta\,z^\beta
\\
0
&
=
\ \ \ \ \ \ 
a_{\epsilon_i}{\sf D}(z_i)
+\cdots\cdots+
a_{n\epsilon_i}{\sf D}(z_i^n)
+
\sum_{\beta\neq\epsilon_i,\dots,n\epsilon_i\atop
1\leqslant\vert\beta\vert\leqslant d}\,
a_\beta\,{\sf D}(z^\beta)
\\
\cdot\cdot
&
\cdots\cdots\cdots\cdots\cdots\cdots\cdots\cdots\cdots
\cdots\cdots\cdots\cdots\cdots\cdots\cdots\cdots\cdot\cdot
\\
0
&
=
\ \ \ \ \ \
a_{\epsilon_i}{\sf D}^n(z_i)
+\cdots+
a_{n\epsilon_i}{\sf D}^n(z_i^n)
+
\sum_{\beta\neq\epsilon_i,\dots,n\epsilon_i\atop
1\leqslant\vert\beta\vert\leqslant d}\,
a_\beta\,{\sf D}^n(z^\beta).
\endaligned\right.
\]
Here in the last $n$ lines, we emphasize a generalized $n \times n$
Wronskian-like matrix, about which the next lemma states that its
determinant is nonzero {\em if and only if} $z_i' \neq 0$, an
assumption we made. For short in the sequel, we shall write $(z_i^k)^{
(\kappa)}$ instead of $D^\kappa ( z_i^k)$, since it is now clear and
unambiguous that primes denote abstract jet variables.

\def\thelemma{\!}\begin{lemma}
For every $i = 1, 2, \dots, n+1$, one has:
\[
\aligned
\left\vert
\begin{array}{cccc}
z_i' & (z_i^2)' & \cdots & (z_i^n)' 
\\
z_i'' & (z_i^2)'' & \cdots & (z_i^n)''
\\
\cdot\cdot & \cdot\cdot & \cdots & \cdot\cdot
\\
z_i^{(n)} & (z_i^2)^{(n)} & \cdots & (z_i^n)^{(n)}
\end{array}
\right\vert
&
=
1!2!\cdots n!\cdot z_i'\,(z_i')^2\cdots(z_i')^n
\\
&
=
1!2!\cdots n!\cdot (z_i')^{\frac{n(n+1)}{2}}.
\endaligned
\]
\end{lemma}

Our appendix is devoted to the proof of this elementary, but
not straightforward, determinantal
identity. As a result, we immediately deduce that $J_{\sf vert}^n (
\mathcal{ X}_0)$ is smooth of pure codimension $(n+1)$ at each one of
its points which lies in $\{ z_i' \neq 0 \}$, for the last $n$
defining equations written above can at first be solved with respect
to $a_{ \epsilon_i}, \dots, a_{ n \epsilon_i}$ thanks to Cramer's
rule, while the first defining equation is trivially solvable with
respect to $a_0$. In other words, in our open set $\{ z_i ' \neq 0\}$,
the vertical affine jet manifold may be represented as a plain 
semi-global {\em graph}:
\[
a_0,a_{\epsilon_i},\dots,a_{n\epsilon_i}
=
\text{\rm certain functions of}\,
\big(
z,z',\dots,z^{(n)},\widetilde{a}_i
\big),
\]
where $\widetilde{ a}_i := \big( a_\beta \big)_{ \beta \neq 0,
\epsilon_i, \dots, n\epsilon_i}^{ 1 \leqslant \vert \beta \vert
\leqslant d}$ gathers all the other coefficients of the universal
hypersurface.

Now, we seek vector fields of the form ${\sf T} = \sum_{ \vert \alpha
\vert \leqslant n}\, {\sf A}_\alpha \, \frac{ \partial }{ \partial
a_\alpha}$ which would be tangent to $J_{ \sf vert}^n ( \mathcal{ X}_0
)$ with the length of the appearing multiindices
being bounded by $n$. For this reason, and because the
equations of $J_{ \sf vert}^n ( \mathcal{ X}_0)$
are {\em linear} with respect to the coefficients $a_\beta$,
when one applies such a derivation ${\sf T}$
to the equations in question, every
monomial $z^\beta$ with $n+1 \leqslant \vert \beta \vert \leqslant d$
disappears automatically, hence we come down 
to solving the following linear system:
\[
\left\{
\aligned
0
&
=
{\sf A}_0
+
{\sf A}_{\epsilon_i}z_i
+\cdots\cdots\cdot+
{\sf A}_{n\epsilon_i}z_i^n
+
\sum_{\alpha\neq\epsilon_i,\dots,n\epsilon_i\atop
1\leqslant\vert\alpha\vert\leqslant n}\,
{\sf A}_\alpha\,z^\alpha
\\
0
&
=
\ \ \ \ \ \ \ \ \
{\sf A}_{\epsilon_i}z_i'
+\cdots\cdot+
{\sf A}_{n\epsilon_i}
(z_i^n)'
+
\sum_{\alpha\neq\epsilon_i,\dots,n\epsilon_i\atop
1\leqslant\vert\alpha\vert\leqslant n}\,
{\sf A}_\alpha\,(z^\alpha)'
\\
\cdot\cdot
&
\cdots\cdots\cdots\cdots\cdots\cdots\cdots\cdots
\cdots\cdots\cdots\cdots\cdots\cdots\cdots\cdots\cdot
\\
0
&
=
\ \ \ \ \ \
{\sf A}_{\epsilon_i}z_i^{(n)}
+\cdots+
{\sf A}_{n\epsilon_i}
(z_i^n)^{(n)}
+
\sum_{\alpha\neq\epsilon_i,\dots,n\epsilon_i\atop
1\leqslant\vert\alpha\vert\leqslant n}\,
{\sf A}_\alpha(z^\alpha)^{(n)},
\endaligned\right.
\]
having the ${\sf A}_\alpha$ as unknowns, where notably, 
$\vert \alpha \vert \leqslant n$ everywhere.

Noticing that the number of directions $\frac{ \partial }{ \partial
a_\alpha}$ equals $\frac{ (n+1+n)!}{ (n+1)! \, n!}$ while the number
of equations above equals $(n+1)$, we may now claim that for every
$\alpha \neq 0, \epsilon_i, \dots, n \epsilon_i$ with
$\vert \alpha \vert \leqslant n$, there are $\frac{
(n+1+n)!}{ (n+1)!  \, n!} - (n+1)$ linearly independent vector fields
of the specific form:
\[
\frac{\partial}{\partial a_\alpha}
-
{\sf B}_{\alpha,0}^i\,
\frac{\partial}{\partial a_0}
-
{\sf B}_{\alpha,1}^i\,
\frac{\partial}{\partial a_{\epsilon_i}}
-\cdots-
{\sf B}_{\alpha,n}^i\,
\frac{\partial}{\partial a_{n\epsilon_i}}
\]
that are tangent to the (semi-global) graph $J_{\sf vert}^n (
\mathcal{ X}_0) \cap \{ z_i' \neq 0\}$, that is to say, the
coefficients of which satisfy the written linear system.  In order to
insure meromorphic prolongation to projective spaces ({\em see}
below), it is convenient to multiply in advance the
basic vector field
$\frac{ \partial }{ \partial a_\alpha}$ of such a kind of
sought vector field
by the Wronskian-like determinant $\Delta ( z_i') := 1! 2! \cdots n!
\cdot (z_i')^{ \frac{ n( n+1)}{ 2}}$ 
that the lemma computed.
In sum, for any $\alpha$ with $1 \leqslant \vert \alpha \vert
\leqslant n$ which is different from $0, \epsilon_i, \dots,
n\epsilon_i$, the vector field:
\[
{\sf T}_\alpha
:=
\Delta(z_i')
\,\frac{\partial}{\partial a_\alpha}
-
{\sf B}_{\alpha,0}^i\,
\frac{\partial}{\partial a_0}
-
{\sf B}_{\alpha,1}^i\,
\frac{\partial}{\partial a_{\epsilon_i}}
-\cdots-
{\sf B}_{\alpha,n}^i\,
\frac{\partial}{\partial a_{n\epsilon_i}}
\]
is tangent to $J_{\sf vert}^n ( \mathcal{ X}_0) \cap
\{ z_i' \neq 0\}$ if and only if its unknown coefficients
${\sf B}_{ \alpha, k}^i$ satisfy the following linear system:
\[
\left\{
\aligned
0
&
=
-{\sf B}_{\alpha,0}^i
-
{\sf B}_{\alpha,1}^i\,z_i
-\cdots\cdot\cdot-
{\sf B}_{\alpha,n}^i\,z_i^n
+
\Delta(z_i')\cdot z^\alpha
\\
0
&
=
\ \ \ \ \ \ \ \ \ \ \
-{\sf B}_{\alpha,1}^i\,z_i'
-\cdots\cdot\cdot-
{\sf B}_{\alpha,n}^i\,(z_i^n)'
+
\Delta(z_i')\cdot (z^\alpha)'
\\
\cdot\cdot
&
\cdots\cdots\cdots\cdots\cdots\cdots\cdots\cdots\cdots
\cdots\cdots\cdots\cdots\cdots\cdot\cdot\cdots\cdots
\\
0
&
=
\ \ \ \ \ \ \ \ \ \ \
-{\sf B}_{\alpha,1}^i\,z_i^{(n)}
-\cdots-
{\sf B}_{\alpha,n}^i\,
(z_i^n)^{(n)}
+
\Delta(z_i')\cdot (z^\alpha)^{(n)}.
\endaligned\right.
\]
A basic application of Cramer's rule now enable us to solve
the last $n$ equations, and afterwards, we may then substitute the
obtained solutions in the first equation:
\[
\left[
\aligned
{\sf B}_{\alpha,k}^i
&
:=
\left\vert
\begin{array}{ccccc}
z_i' & \cdots & (z^\alpha)' & \cdots & (z_i^n)'
\\
\vdots & \ddots & \vdots & \ddots & \vdots
\\
z_i^{(n)} & \cdots & (z^\alpha)^{(n)} & \cdots & (z_i^n)^{(n)}
\end{array}
\right\vert
\ \ \ \ \ \ \ \ \ \ \ \ \
{\scriptstyle{(\kappa\text{\rm -th column},\,\,\,
1\,\leqslant\,k\,\leqslant\,n)}}
\\
{\sf B}_{\alpha,0}^i
&
:=
-{\sf B}_{\alpha,1}^i\,z_i
-\cdots-
{\sf B}_{\alpha,n}^i\,z_i^n
+
\Delta(z_i')\cdot z^\alpha.
\endaligned\right.
\]
Clearly, the so obtained vector fields ${\sf T}_\alpha$ with $\alpha
\neq 0, \epsilon_i, \dots, n\epsilon_i$ are linearly independent at
every point of $J_{ \sf vert}^n ( \mathcal{ X}_0) \cap \{ z_i' = 0
\}$.

\subsection*{ Meromorphic prolongation and computation of
pole orders}
Recall that any polynomial $P ( t_0, \dots, 
t_{ \upsilon-1}, t_{ \upsilon +1}, \dots, t_{n+1})$ of degree $e
\geqslant 1$ on an affine $\C^{ n + 1} \subset \P^{ n + 1}$, when
viewed as a meromorphic map $\P^{ n + 1} \to \P^1$, has pole order
equal to $e$, for a change of standard affine chart:
\[
t_0
=
{\textstyle{\frac{1}{z_\upsilon}}},
\dots,
t_{\upsilon-1}
=
{\textstyle{\frac{z_{\upsilon-1}}{z_\upsilon}}},\,\,
t_{\upsilon+1}
=
{\textstyle{\frac{z_{\upsilon+1}}{z_\upsilon}}},
\dots,
t_{n+1}
=
{\textstyle{\frac{z_{n+1}}{z_\upsilon}}},
\]
transfers $P$ to $P \big( \frac{ 1}{ z_\upsilon}, \dots, \frac{ z_{
\upsilon-1}}{ z_\upsilon}, \frac{ z_{ \upsilon +1}}{ z_\upsilon},
\dots, \frac{ z_{ n+1}}{ z_\upsilon} \big)$. Through such an inversion
map, the first-order jets, second-order jets, 
{\em etc}, are transferred to:
\[
\aligned
{\textstyle{\frac{z_i'}{z_\upsilon}}}
-
{\textstyle{\frac{z_iz_\upsilon'}{z_\upsilon^2}}},
\ \ \ \ \ \ \ \ \ \ \ \
{\textstyle{\frac{z_i''}{z_\upsilon}}}
-
2\,{\textstyle{\frac{z_i'z_\upsilon'}{z_\upsilon^2}}}
-
{\textstyle{\frac{z_iz_\upsilon''}{z_\upsilon^2}}}
+
2\,{\textstyle{\frac{z_iz_\upsilon'z_\upsilon'}{z_\upsilon^3}}},
\ \ \ \ \ \
\text{\em etc.},
\endaligned
\]
hence by just looking at the maximal power of $z_\upsilon$ at the
denominator, one easily observes by induction that:
\[
\text{\sf Pole-order}
\left[
z^\alpha\,
(z')^{\alpha^1}\cdots
(z^{(n)})^{\alpha^n}
\right]
=
\vert\alpha\vert
+
\vert\alpha^1\vert
+\cdots+
\vert\alpha^n\vert
+
n,
\]
and furthermore, one differentiation of such a monomial
increases its pole order by just one unit. Now, we claim that:
\[
\left[
\aligned
\text{\sf Pole-order}
\big[\Delta(z_i')\big]
&
=
n^2+n,
\\
\text{\sf Pole-order}
\big[
{\sf B}_{\alpha,k}^i
\big]
&
=
\vert\alpha\vert
+
n^2+n-k
\\
\text{\sf Pole-order}
\big[
{\sf B}_{\alpha,0}^i
\big]
&
=
\vert\alpha\vert+n^2+n,
\endaligned\right.
\]
so that the highest pole order occurs to be the coefficient 
${\sf B}_{\alpha,0}^i$
of $\frac{
\partial }{ \partial a_0}$ in each ${\sf T}_\alpha$.

Indeed, replacing the entries of the determinant $\Delta ( z_i')$
plainly by the nonnegative integers which indicate the pole orders, we
may write symbolically:
\[
\text{\sf Pole-order}
\big[\Delta(z_i')\big]
=
\text{\sf Pole-order of}\,\,
\left\vert
\begin{array}{ccccc}
2 & 3 & 4 & \cdots & n+1
\\
3 & 4 & 5 & \cdots & n+2
\\
4 & 5 & 6 & \cdots & n+3
\\
\cdot\cdot & \cdot\cdot & \cdot\cdot & \cdots & \cdot\cdot
\\
n+1 & n+2 & n+3 & \cdots & 2n
\end{array}
\right\vert.
\]
When one expands the determinant as a sum of monomials with $\pm$ signs,
pole orders are just added, symbolically speaking.  Then one easily
convinces oneself that {\em each one} of the obtained monomials has
the {\em same} pole order, hence it suffices
to compute the pole order of the
monomial of the main diagonal, which 
is equal to: $2 + 4 + 6 + \cdots + 2n = n(n+1)$.

Next, ${\sf B}_{\alpha, k}^i$ is obtained from $\Delta ( z_i')$
by replacing the $k$-th column of $\Delta ( z_i')$ by the 
new column of
pole order entries $\vert \alpha \vert + 1$, $\vert \alpha \vert + 2$,
\dots, $\vert \alpha \vert + n$. The pole order $\vert \alpha \vert$
being ``factorizable'', we get:
\[
\text{\sf Pole-order}
\big[{\sf B}_{\alpha,k}^i\big]
=
\vert\alpha\vert
+
\text{\sf Pole-order of}\,\,
\left\vert
\begin{array}{cccccc}
2 & 3 & \cdots & 1 & \cdots & n+1
\\
3 & 4 & \cdots & 2 & \cdots & n+2
\\
\cdot\cdot & \cdot\cdot & \cdots & \cdot\cdot &
\cdots & \cdot\cdot
\\
n+1 & n+2 & \cdots & n & \cdots & 2n 
\end{array}
\right\vert,
\]
where the central-looking column is the $k$-th, the only which
differs from $\Delta ( z_i')$. Again, one easily convinces oneself
that the pole order of {\em each one} of the monomials obtained after
expansion is the {\em same}, so that by looking again
at the main diagonal:
\[
\aligned
\text{\sf Pole-order}
\big[{\sf B}_{\alpha,k}^i\big]
&
=
\vert\alpha\vert
+
2+\cdots+
2(k-1)+
k+
2(k+1)
+\cdots+2n
\\
&
=
\vert\alpha\vert
+
n(n+1)
-k.
\endaligned
\]

Finally, coming back to the definition of ${\sf B}_{ \alpha, 0}^i$,
one then immediately sees that each term in ${\sf B}_{ \alpha, 0}^i$
has the same pole order, equal to $\vert \alpha \vert + n^2 + n$.

In conclusion, the maximal pole order is reached by any ${\sf B}_{
\alpha, 0}^i$ with $\vert \alpha \vert = n$, and is equal to $n^2 +
2n$, as it appears in the statement of the main theorem.  Clearly, the
poles are compensated by the twisting $( \bullet) \otimes \mathcal{
O}_{ \P^{ n+1}} ( n^2 + 2n)$. 
Notice that
the coefficients of the constructed
${\sf T}_\alpha$'s depend only on the jet variables $(z, z', \dots,
z^{ (n)})$, absolutely not on the coefficients $a_\beta$.

The other vector fields that we will construct in the remainder of
the article, so as to complete a true framing, will all have pole
order in the $z$-direction smaller than $n^2 + 2n$, and will all
have pole order at most $1$ in the $a$-direction.

\subsection*{ Modifications needed for the second theorem}
By assumption, at least one (classical) Wronskian $\det \big( z_i^{ (
\lambda_j)} \big)_{ 1 \leqslant i \leqslant n+1}^{ 1 \leqslant j
\leqslant n}$ does not vanish, where, without loss of generality, we
may assume that $1 \leqslant \lambda_1 < \cdots < \lambda_n \leqslant
n+1$. To fix ideas, we shall work in the open set $\big\{ \det (
z_i^{(j)})_{ 1 \leqslant i \leqslant n}^{ 1 \leqslant j \leqslant n}
\neq 0 \big\}$ and we shall denote: ${\sf W} := \det ( z_i^{(j)})_{ 1
\leqslant i \leqslant n}^{ 1 \leqslant j \leqslant n}$.  Then again,
the variety of vertical jets is checked to be, in the open set $\big\{
\det ( z_i^{(j)})_{ 1 \leqslant i \leqslant n}^{ 1 \leqslant j
\leqslant n} \neq 0 \big\}$, a semi-global graph of equations:
\[
\left\{
\aligned
0
&
=
a_0+a_{\epsilon_1}z_1
+\cdots\cdot\cdot+
a_{\epsilon_n}z_n
+
\sum_{\beta\neq\epsilon_1,\dots,\epsilon_n\atop
1\leqslant\vert\beta\vert\leqslant d}\,
a_\beta\,z^\beta
\\
0
&
=
\ \ \ \ \ \ \ \ \
a_{\epsilon_1}z_1'
+\cdots\cdot\cdot+
a_{\epsilon_n}z_n'
+ 
\sum_{\beta\neq\epsilon_1,\dots,\epsilon_n\atop
1\leqslant\vert\beta\vert\leqslant d}\,
a_\beta\,(z^\beta)'
\\
\cdot\cdot
&
\cdots\cdots\cdots\cdots\cdots\cdots\cdots\cdots\cdots
\cdots\cdots\cdots\cdots\cdots\cdots
\\
0
&
=
\ \ \ \ \ \ \
a_{\epsilon_1}z_1^{(n)}
+\cdots+
a_{\epsilon_n}z_n^{(n)}
+
\sum_{\beta\neq\epsilon_1,\dots,\epsilon_n\atop
1\leqslant\vert\beta\vert\leqslant d}\,
a_\beta\,(z^\beta)^{(n)},
\endaligned\right.
\]
having transversal coordinates $(a_0, a_{ \epsilon_1}, \dots, a_{
\epsilon_n})$, the ones that are then 
clearly solvable here. Thus if,
similarly as in the previous paragraphs, one seeks tangent vector
fields of the specific form:
\[
{\sf T}_\alpha
:=
{\sf W}\,\frac{\partial}{\partial a_\alpha}
-
{\sf B}_{\alpha,0}\,\frac{\partial}{\partial a_0}
-
{\sf B}_{\alpha,1}\,\frac{\partial}{\partial a_{\epsilon_1}}
-\cdots-
{\sf B}_{\alpha,n}\,\frac{\partial}{\partial a_{\epsilon_n}},
\]
for any $\alpha \in \N^{ n+1}$ with $\vert \alpha \vert \leqslant n$
and with $\alpha \neq 0, \epsilon_1, \dots, \epsilon_n$, then the
system we now have to solve becomes:
\[
\left\{
\aligned
0
&
=
-{\sf B}_{\alpha,0}
-
{\sf B}_{\alpha,1}\,z_1
-\cdots\cdots-
{\sf B}_{\alpha,n}\,z_n
+
{\sf W}\cdot z^\alpha
\\
0
&
=
\ \ \ \ \ \ \ \ \ \ \
-{\sf B}_{\alpha,1}\,z_1'
-\cdots\cdots-
{\sf B}_{\alpha,n}\,z_n'
+
{\sf W}\cdot(z^\alpha)'
\\
\cdot\cdot
&
\cdots\cdots\cdots\cdots\cdots\cdots\cdots\cdots\cdots
\cdots\cdots\cdots\cdots\cdots\cdots
\\
0
&
=
\ \ \ \ \ \ \ \ \ \
-{\sf B}_{\alpha,1}\,z_1^{(n)}
-\cdots\cdot-
{\sf B}_{\alpha,n}\,z_n^{(n)}
+
{\sf W}\cdot(z^\alpha)^{(n)}. 
\endaligned\right.
\]
The unique solution is then again yielded by Cramer's rule:
\[
\left[
\aligned
{\sf B}_{\alpha,k}
&
:=
\left\vert
\begin{array}{ccccc}
z_1' & \cdots & (z^\alpha)' & \cdots & z_n'
\\
\vdots & \ddots & \vdots & \ddots & \vdots
\\
z_1^{(n)} & \cdots & (z^\alpha)^{(n)} & \cdots &
z_n^{(n)}
\end{array}
\right\vert
\ \ \ \ \ \ \ \ \ \ \ \ \
{\scriptstyle{(k\text{\rm -th column},\,\,\,
1\,\leqslant\,k\,\leqslant\,n)}}
\\
{\sf B}_{\alpha,0}
&
:=
-{\sf B}_{\alpha,1}\,z_1
-\cdots-
{\sf B}_{\alpha,n}\,z_n
+
{\sf W}\cdot z^\alpha.
\endaligned\right.
\]
One may now verify that:
\[
\aligned
\text{\sf Pole order}
\big[
{\sf W}
\big]
&
=
{\textstyle{\frac{(n+1)(n+2)}{2}}}
\\
\text{\sf Pole order}
\big[
{\sf B}_{\alpha,k}
\big]
&
=
2+\cdots+(n+1)
-(k+1)+\vert\alpha\vert+k
\\
&
=
{\textstyle{\frac{n^2+3n-2}{2}}}
+\vert\alpha\vert
\\
\text{\sf Pole order}
\big[
{\sf B}_{\alpha,0}
\big]
&
=
{\textstyle{\frac{n^2+3n}{2}}}
+\vert\alpha\vert,
\endaligned
\]
so that the maximal pole order is reached by ${\sf B}_{ \alpha, 0}$
for any multiindex $\alpha$ with $\vert \alpha \vert = n$, and is
equal to $\frac{ n^2 + 5n }{ 2}$, as this appears in the second
theorem.

The other vector fields that we will construct in the sequel will
complete a generating set both for the first and for the second
theorem and will have lower pole order in the $z$-direction.

\subsection*{ Higher lengths}
At present, we construct globally defined tangent vector fields
which span the remaining directions $\bigoplus_{ n+1 \leqslant \vert
\alpha \vert \leqslant d \atop \alpha_1 < d}\, \C\cdot \frac{ \partial
}{\partial a_\alpha}$ in the space of coefficients $a_\alpha$.  For an
arbitrary multiindex $\ell = (\ell_1, \ell_2, \dots, \ell_{ n+1} ) \in
\N^{ n + 1}$ of length:
\[
n+1
=
\ell_1+\ell_2+\cdots+\ell_{n+1},
\]
we introduce the following family of vector fields living only in the
space of $a$-variables:
\[
{\sf T}_\alpha^{\ell_1,\ell_2,\dots,\ell_{n+1}}
=
{\sf T}_\alpha^\ell
=
\sum_{\ell'+\ell''=\ell\atop\ell',\,\ell''\in\N^{n+1}}\,
(-1)^{\vert\ell''\vert}\,
\frac{\ell!}{\ell'!\,\,\ell''!}\,
z^{\ell''}\,
\frac{\partial}{\partial a_{\alpha-\ell''}},
\]
where the indices $\alpha$ are all possible indices satisfying
$\alpha_1 \geqslant \ell_1$, \dots, $\alpha_{ n+1} \geqslant \ell_{
n+1}$, $\vert \alpha \vert \leqslant d$ and $\alpha_1 < d$, and where
the sum abbreviates $\sum_{ \ell_1 ' + \ell_1 '' = \ell_1} \cdots\,
\sum_{ \ell_{ n+1}' + \ell_{ n+1}'' = \ell_{ n+1}}$.  For instance,
for $n+1 = 4$ and with the special choice $\ell_1 = \ell_2 = 2$
(whence necessarily $\ell_3 = \ell_4 = 0$), we get the following
family of vector fields defined for all $\alpha$ with $\alpha_1
\geqslant 2$, $\alpha_2 \geqslant 2$, $\alpha_3 \geqslant 0$,
$\alpha_4 \geqslant 0$ and $\vert \alpha \vert \leqslant d$, $\alpha_1
< d$ (compare~\cite{ rou2007a}, p.~373):
\[
\aligned
{\sf T}_\alpha^{2,2,0,0}
&
=
\frac{\partial}{\partial a_\alpha}
-
2z_1\,\frac{\partial}{\partial a_{\alpha-\epsilon_1}}
-
2z_2\,\frac{\partial}{\partial a_{\alpha-\epsilon_2}}
+
\\
&
\ \ \ \ \
+
z_1^2\,\frac{\partial}{\partial a_{\alpha-2\epsilon_1}}
+
4z_1z_2\,\frac{\partial}{\partial a_{\alpha-\epsilon_1-\epsilon_2}}
+
z_2^2\,\frac{\partial}{\partial a_{\alpha-2\epsilon_2}}
-
\\
&
\ \ \ \ \
-
2z_1^2z_2\,\frac{\partial}{\partial a_{\alpha-2\epsilon_1-\epsilon_2}}
-
2z_1z_2^2\,\frac{\partial}{\partial a_{\alpha-\epsilon_1-2\epsilon_2}}
+
z_1^2z_2^2\,\frac{\partial}{\partial a_{\alpha-2\epsilon_1-2\epsilon_2}}.
\endaligned
\]
After a moment's reflection, one may convince oneself that as $\ell$
with $\vert \ell \vert = n + 1$ runs and as $\alpha$ with $\alpha_i
\geqslant \ell_i$ runs, the ${\sf T}_\alpha^\ell$ together with the
vector fields of the previous paragraph do span $\bigoplus_{ n +1
\leqslant \vert \alpha \vert \leqslant d \atop \alpha_1 < d}\, \C
\cdot \frac{\partial }{\partial a_\alpha}$; there are in fact
redundancies among the {\em triangular} system defined by the ${\sf
T}_\alpha^\ell$, whenever one has $\alpha \geqslant \ell_1$ and
$\alpha \geqslant \ell_2$ for two distinct $\ell^1$, $\ell^2$ with
$\vert \ell^1 \vert = \vert \ell^2 \vert = n + 1$.

\def\thelemma{\!\!}\begin{lemma}
For every nonnegative integer $e \leqslant n$ and for arbitrary
indices $j_1, \dots, j_e$ with $1 \leqslant j_i \leqslant n+1$, one
has:
\[
0
\equiv
{\sf T}_\alpha^\ell
\bigg(
\sum_{\beta\in\N^{n+1}\atop\vert\beta\vert\leqslant d,\,a_{d0\cdots 0}
=1}\,
a_\beta\,
\frac{\partial^e(z^\beta)}{\partial z_{j_1}\cdots
\partial z_{j_e}}
\bigg),
\]
and as a result, ${\sf T}_\alpha^\ell$ identically annihilates all the
defining equations of $J_{\sf vert}^n ( \mathcal{ X}_0)$, hence is
tangent to $J_{\sf vert}^n ( \mathcal{ X}_0)$.
\end{lemma}

\proof
Let $w_1, w_2, \dots, w_{ n+1}$ be auxiliary complex variables.  For
every derivation order $e \leqslant n$ strictly less than the
vanishing order $\sum\, \ell_i = n +1$, we trivially have:
\[
0
\equiv
\frac{\partial}{\partial z_{j_1}}\ldots
\frac{\partial}{\partial z_{j_e}}
\Big(
\left[w_1-z_1\right]^{\ell_1}
\left[w_2-z_2\right]^{\ell_2}
\cdots
\left[w_{n+1}-z_{n+1}\right]^{\ell_{n+1}}
\Big)\Big\vert_{w=z}.
\]
In other words, by expanding $\left[ w - z \right]^\ell = \sum_{\ell '
+ \ell '' = \ell}\, (-1)^{ \vert \ell '' \vert} \, \frac{ \ell! }{\ell
' ! \, \, \ell''!} \, w^{ \ell'} \, z^{ \ell''}$ thanks to the
multinomial formula, by letting the derivation $\partial ^e ( \bullet
) \big/ \partial z_{ j_1} \cdots \partial z_{ j_e}$ act on this
expansion, by setting $w = z$, and finally, by multiplying the result
obtained by $z^{ \alpha - \ell}$, we get the useful identities:
\[
\boxed{
0
\equiv
\sum_{\ell'+\ell''=\ell\atop\ell',\,\ell'\in\N^{n+1}}
(-1)^{\vert\ell''\vert}\,
\frac{\ell!}{\ell'!\,\,\ell''!}\,
z^{\alpha-\ell''}\cdot
\frac{\partial^e\big(z^{\ell''}\big)}{\partial z_{j_1}\cdots
\partial z_{j_e}}
}\,.
\]
On the other hand, by letting the derivation ${\sf T}_\alpha^\ell$ act
as it should, the identities of the lemma that we have to check may be
written:
\[
\aligned
0
&
\overset{?}{\equiv}
\sum_{\ell'+\ell''=\ell\atop\ell',\,\ell'\in\N^{n+1}}\,
(-1)^{\vert\ell''\vert}\,
\frac{\ell!}{\ell'!\,\,\ell''!}\,
\frac{\partial}{\partial a_{\alpha-\ell''}}\,
\bigg(
\sum_{\beta\in\N^{n+1}\atop
\vert\beta\vert\leqslant d,\,a_{d0\cdots0}=1}\,
a_\beta\,
\frac{\partial^e(z^\beta)}{\partial z_{j_1}\cdots\partial z_{j_e}}
\bigg)\cdot
z^{\ell''}
\\
&
=
\sum_{\ell'+\ell''=\ell\atop\ell',\,\ell'\in\N^{n+1}}\,
(-1)^{\vert\ell''\vert}\,
\frac{\ell!}{\ell'!\,\,\ell''!}\,
\frac{\partial^e\big(z^{\alpha-\ell''}\big)}{
\partial z_{j_1}\cdots\partial z_{j_e}}\cdot
z^{\ell''}.
\endaligned
\]
Compared to the boxed, known identities, the derivation is now
switched to the other monomial. Generally, we claim that for every $e
= 0, 1, \dots, n$ and for every decomposition $e = e_1 + (e - e_1)$
with $0 \leqslant e_1 \leqslant e$, the expression:
\[
\small
\aligned
\big(
j_1,\dots,j_{e_1}\big\vert
j_{e_1+1},\dots,j_e\big)
:=
\sum_{\ell'+\ell''=\ell\atop\ell',\,\ell'\in\N^{n+1}}\,
(-1)^{\vert\ell''\vert}\,
\frac{\ell!}{\ell'!\,\,\ell''!}\,
\frac{\partial^{e_1}\big(z^{\alpha-\ell''}\big)}{
\partial z_{j_1}\cdots\partial z_{j_{e_1}}}\cdot
\frac{\partial^{e-e_1}\big(z^{\ell''}\big)}{
\partial z_{j_{e_1+1}}\cdots\partial z_{j_e}}
\endaligned
\]
vanishes identically, for all indices $j_1, \dots, j_e = 1, 2, \dots
n+1$.  We know that this assertion is true when $e_1 = 0$ for all $e =
0, 1, \dots, n$ and the lemma corresponds to $e - e_1 = 0$ for all $e_1 =
0, 1, \dots, n$.

For $e = 0$, the assertion is thus known.  Suppose it to be true at
level $e$. Reasoning by induction, we then assume that:
\[
0
\equiv
\big(
j_1,\dots,j_{e_1}\big\vert
j_{e_1+1},\dots,j_e
\big),
\]
for all $e_1 = 0, 1, \dots, e$ and all possible $j_i$, If $e + 1$ is
still $\leqslant n$, we differentiate all these identities with respect
to $z_k$ using Leibniz' rule and we organize the resulting equations
as a convenient array:
\[
\aligned
0
&
\equiv
\big(j_1,\dots,j_e,k\big\vert\emptyset\big)
+
\big(j_1,\dots,j_e\big\vert k\big)
\\
0
&
\equiv
\big(j_1,\dots,j_{e-1},k\big\vert j_e\big)
+
\big(j_1,\dots,j_{e-1}\big\vert j_e,k\big)
\\
&
\cdots\cdots\cdots\cdots\cdots\cdots\cdots\cdots\cdots\cdots\cdots
\cdots\cdots
\\
0
&
\equiv
\big(j_1,k\big\vert j_2,\dots,j_e\big)
+
\big(j_1\big\vert j_2,\dots,j_e,k\big)
\\
0
&
\equiv
\big(k\big\vert j_1,\dots,j_e\big)
+
\underline{\big(\emptyset\big\vert j_1,\dots,j_e,k\big)}_{\circ}.
\endaligned
\]
We have underlined the last term, known to vanish. Then the first term
of the last line vanishes, for all indices $k, j_1, \dots, j_e = 1, 2,
\dots, n+1$.  So the second term of the penultimate vanishes, {\em
etc.}, and hence the very first term $\big( j_1, \dots, j_e, k
\big\vert \emptyset \big)$ does vanish identically, as desired.
\endproof

\section*{\S4.~Second package of jet, coordinate vector fields}
\label{Section-4}

\subsection*{ Spanning the $\frac{ \partial }{\partial
z_i}$-directions}
To complete the framing, let us at first span all
the $\frac{ \partial }{\partial z_i}$ directions.
By convention, $a_{ d0\cdots 0} = 1$.

\def\thelemma{\!}\begin{lemma}
For $i = 1, 2, \dots, n+1$, the vector fields:
\[
{\sf T}_i
:=
\frac{\partial}{\partial z_i}
-
\sum_{\vert\alpha\vert\leqslant d-1}\,
a_{\alpha+\epsilon_i}(\alpha_i+1)\,
\frac{\partial}{\partial a_\alpha}
\]
are all tangent to $J_{\sf vert}^n ( \mathcal{ X}_0)$.
\end{lemma}

\proof
Appying the derivation ${\sf T}_i$ to the first equation $0 =
\sum_\alpha \, a_\alpha z^\alpha$ of $J_{\sf vert}^n ( \mathcal{
X}_0)$, we indeed get an identically vanishing result:
\[
\sum_{\vert\alpha\vert\leqslant d}\,
a_\alpha\,\frac{\partial(z^\alpha)}{\partial z_i}
-
\sum_{\vert\alpha\vert\leqslant d-1}\,
a_{\alpha+\epsilon_i}(\alpha_i+1)\,z^\alpha
\equiv
0.
\]
Since the ${\sf T}_i$ commute with the total differentiation operator
${\sf D}$, it then follows immediately that ${\sf T}_i$ annihilates
all the other defining equations:
\[
0
\equiv
{\sf T}_i
\Big(
{\sf D}\,
\sum_\alpha\,a_\alpha 
z^\alpha
\Big)
\equiv
\cdots
\equiv
{\sf T}_i
\Big(
{\sf D}^n\,
\sum_\alpha\,a_\alpha 
z^\alpha
\Big),
\]
and this yields the tangency property claimed.
\endproof

\subsection*{Spanning the 
$\partial \big/ \partial z_j^{( \lambda )}$ directions} For the last
family of vector fields, we transfer to general $\kappa = n \geqslant
2$ the approach of \cite{ pau2008} known for $\kappa = n = 2$ and also
for $\kappa = n =3$ \cite{ rou2007a}, with few differences.
 
Let $\Lambda = \big( \Lambda_k^l \big)_{ 1 \leqslant k \leqslant
n+1}^{ 1 \leqslant l \leqslant n + 1}$ be a matrix in ${\sf GL} ( n +
1, \C)$.  To span the only remaining directions $\partial \big/ \partial
z_j^{( \lambda )}$, one seeks meromorphic vector fields tangent
to $J_{\sf vert}^n ( \mathcal{ X}_0) \big\backslash
\Sigma_0$ that are of the special
form:
\[
\small
\aligned
{\sf T}_\Lambda
&
:=
\sum_{k=1}^{n+1}\,
\Big(
{\textstyle{\sum}_{l=1}^{n+1}}\,
\Lambda_k^l\,z_l'
\Big)\,
\frac{\partial}{\partial z_k'}
+\cdots+
\sum_{k=1}^{n+1}\,
\Big(
{\textstyle{\sum}_{l=1}^{n+1}}\,
\Lambda_k^l\,z_l^{(n)}
\Big)\,
\frac{\partial}{\partial z_k^{(n)}}
+
\\
&
\ \ \ \ \ \ \ \ 
+
\sum_{\vert\alpha\vert\leqslant d\atop\alpha_1<d}\,
{\sf A}_\alpha\big(z,a,\Lambda)\,
\frac{\partial}{\partial a_\alpha},
\endaligned
\]
where, for various jet orders $\lambda$'s, the coefficients ${\sf
Z}_k^{ ( \lambda)}$ of the $\frac{ \partial }{ \partial z_k^{
(\lambda)}}$, $k = 1, \dots, n+1$, are defined a priori to be obtained
by multiplying the jet matrix $\big( z_j^{ ( \lambda)} \big)_{ 1
\leqslant j \leqslant n+1}^{ 1 \leqslant \lambda \leqslant n}$ by such
a matrix $\Lambda$:
\[
\small
\aligned
\left(
\begin{array}{cccc}
\Lambda_1^1 & \cdots & \Lambda_1^n & \Lambda_1^{n+1}
\\
\Lambda_2^1 & \cdots & \Lambda_2^n & \Lambda_2^{n+1}
\\
\vdots & \ddots & \vdots & \vdots
\\
\Lambda_{n+1}^1 & \cdots & \Lambda_{n+1}^n & \Lambda_{n+1}^{n+1}
\end{array}
\right)
\,
\left(
\begin{array}{ccc}
z_1' & \cdots & z_1^{(n)}
\\
z_2' & \cdots & z_2^{(n)}
\\
\vdots & \ddots & \vdots
\\
z_{n+1}' & \cdots & z_{n+1}^{(n)}
\end{array}
\right)
=
\left(
\begin{array}{ccc}
{\sf Z}_1' & \cdots & {\sf Z}_1^{(n)}
\\
{\sf Z}_2' & \cdots & {\sf Z}_2^{(n)}
\\
\vdots & \ddots & \vdots
\\
{\sf Z}_{n+1}' & \cdots & {\sf Z}_{n+1}^{(n)}
\end{array}
\right),
\endaligned
\]
and where the coefficients ${\sf A}_\alpha \big( z, a, \Lambda \big)$,
to be computed shortly, should insure that ${\sf T}_\Lambda$ is
effectively tangent to $J_{\sf vert}^n ( \mathcal{ X}_0 ) \big
\backslash \Sigma_0$.

In fact, by plainly inspecting ranks of the matrix multiplication
above, one easily sees that, at every point of our basic open set
where at least one $n \times n$ (sub)Wronskian of the jet matrix
$\big( z_j^{ ( \lambda)} \big)$ does not vanish, one has for $\Lambda$
varying without restriction in ${\sf GL} ( n+1, \C)$:
\[
\text{\rm Span}_\Lambda
\bigg(
\Lambda\,z'\,\frac{\partial}{\partial z'}
+\cdots+
\Lambda\,z^{(n)}\,\frac{\partial}{\partial z^{(n)}}
\bigg)
=
\bigoplus_{1\leqslant k\leqslant n+1}\!\!
\C\,\frac{\partial}{\partial z_k'}
\,\,\,\cdots\,\,
\bigoplus_{1\leqslant k\leqslant n+1}\!\!
\C\,\frac{\partial}{\partial z_k^{(n)}}.
\]  
The following proposition will therefore complete the proof of the
theorem.

\def\theproposition{\!}\begin{proposition}
There exist coefficients ${\sf A}_\alpha$ for $\partial \big/ \partial
a_\alpha$ with $\vert \alpha \vert \leqslant d$, $\alpha_1 <d$, which
are polynomials in $z$ of degree at most $n$:
\[
{\sf A}_\alpha
\big(z,a,\Lambda\big)
=
\sum_{\vert\beta\vert\leqslant n}\,
\mathcal{L}_\alpha^\beta
\big(a,\Lambda\big)\,
z^\beta
\]
with coefficients $\mathcal{ L}_\alpha^\beta \big( a, \Lambda \big)$
being bilinear in the variables $\big( a_\gamma, \, \Lambda_k^l \big)$
such that ${\sf T}_\Lambda$ is tangent to $J_{\sf vert}^n ( \mathcal{
X}_0) \big \backslash \Sigma_0$.
\end{proposition}

\proof
While writing down, say, the first two tangency equations, namely when
applying the derivative ${\sf T}_\Lambda$ to the first two of the five
big equations written at the beginning, one gets equations:
\def\theequation{${\bf 0}$}\begin{equation}
0
=
\sum_{\vert\alpha\vert\leqslant d
\atop\alpha_1<d}\,
{\sf A}_\alpha\cdot z^\alpha
\end{equation}
\def\theequation{${\bf 1}_{j_1}$}\begin{equation}
0
=
\sum_{\vert\alpha\vert\leqslant d\atop 
\alpha_1<d}\,
{\sf A}_\alpha\cdot 
\frac{\partial(z^\alpha)}{\partial z_{j_1}}
+
\sum_{\vert\alpha\vert\leqslant d\atop a_{d0\cdots 0}=1}\!\!
a_\alpha\,\,
\sum_{l=1}^{n+1}\,
\frac{\partial(z^\alpha)}{\partial z_l}\,\Lambda_l^{j_1},
\end{equation}
for which one is allowed to equate to zero the coefficient of each
$z_{ j_1}'$, because the sought ${\sf A}_\alpha$ should be independent
of $z_{ j_1}', z_{ j_2}'', \dots, z_{ j_n}^{ (n)}$.

Next, when applying ${\sf T}_\Lambda$ to the third defining equation
of $J_{\sf vert}^n ( \mathcal{ X}_0)$, one sees thanks to $({\bf 1}
)_{ j_1}$ that the coefficient of each $z_{ j_1}''$ then automatically
vanishes\footnote{\, This simplification trick justifies a posteriori,
{\em cf.}~\cite{ pau2008}, the {\em ad hoc}-looking assumption that
the same matrix $\Lambda$ appears in each jet vector field coefficient
${\sf Z}^{ ( \lambda)} = \Lambda \cdot z^{ ( \lambda)}$.}, hence we
are left with just equating to zero the coefficients of the monomials
$z_{ j_1}' z_{ j_2}'$, namely:
\def\theequation{${\bf 2}_{j_1j_2}$}\begin{equation}
0
=
\sum_{\vert\alpha\vert\leqslant d\atop
\alpha_1<d}\,
{\sf A}_\alpha\cdot 
\frac{\partial^2(z^\alpha)}{\partial z_{j_1}
\partial z_{j_2}}
+
\sum_{\vert\alpha\vert\leqslant d\atop a_{d0\cdots 0}=1}\!\!
a_\alpha\,\,
\sum_{l=1}^{n+1}
\bigg(
\frac{\partial^2(z^\alpha)}{\partial z_l\partial z_{j_2}}\,
\Lambda_l^{j_1}
+
\frac{\partial^2(z^\alpha)}{\partial z_{j_1}\partial z_l}\,
\Lambda_l^{j_2}
\bigg).
\end{equation}
By induction, such a simplification is easily seen to generalize and
thus, the $(e + 1)$-th condition of tangency, after taking account of
the successive cancellations, is obtained by just looking at how
${\sf T}_\Lambda$ acts on the jet monomial $z_{ j_1}' \cdots z_{ j_e}
'$, and the result then consists in the family of equations:
\def\theequation{${\bf e}_{j_1\cdots j_e}$}\begin{equation}
\small
\aligned
0
&
=
\sum_{\vert\alpha\vert\leqslant d
\atop\alpha_1<d}\,
{\sf A}_\alpha\cdot 
\frac{\partial^e(z^\alpha)}{\partial z_{j_1}\cdots 
\partial z_{j_e}}
+
\sum_{\vert\alpha\vert\leqslant d\atop a_{d0\cdots 0}=1}\!\!
a_\alpha\,\,
\sum_{l=1}^{n+1}
\bigg(
\\
&
\ \ \ \ \ \ \
\bigg(
\frac{\partial^e(z^\alpha)}{\partial z_l\partial z_{j_2}
\cdots\partial z_{j_e}}\,\Lambda_l^{j_1}
+
\frac{\partial^e(z^\alpha)}{\partial z_{j_1}\partial z_l
\cdots\partial z_{j_e}}\,\Lambda_l^{j_2}
+\cdots+
\frac{\partial(z^\alpha)}{\partial z_{j_1}\cdots 
\partial z_{j_{e-1}}\partial z_l}\,
\Lambda_l^{j_e}
\bigg),
\endaligned
\end{equation}
where $j_1, \dots, j_e = 1, \dots, n+1$ are arbitrary.

The equations for the unknows $\mathcal{ L}_\alpha^\beta$ shall then
be obtained by identifying the coefficients of the monomials $z^\rho$
in the above equations $( {\bf 0})$, $( {\bf 1}_{ j_1})$, \dots, $(
{\bf n}_{j_1\cdots j_n})$.

At first, we observe that since the degrees in $z$ of the second terms
of $( {\bf 1}_{ j_1})$, $({\bf 2}_{ j_1 j_2})$, {\em etc.}  are at
most $d-1$, $d-2$, {\rm etc.}, we can, without loss of generality,
suppose that the $\mathcal{ L}_\alpha^\beta$ are zero for $\vert
\alpha \vert + \vert \beta \vert \geqslant d +1$, as it is written in
the proposition.  Next ({\rm cf.}~\cite{ pau2008}), using the equation
of $\mathcal{ X}_0$, we may replace the occurence of $z_1^d$ in the
equation ({\bf 0}) by $- \sum_{ \vert\alpha \vert \leqslant d, \,
\alpha_1 < d}\, a_\alpha \, z^\alpha$, so that the degree in the $z_1$
variable is at most $d-1$ (as in~\cite{ pau2008}, this will insure
that the linear systems we have to solve are not overdetermined, and
Cramer's basic rule will apply).

Now, the coefficient of each monomial $z^\rho$ in the equation $({\bf
0})$ should vanish:
\def\theequation{${\bf 0}_\rho$}\begin{equation}
0
=
\sum_{\alpha+\beta=\rho}\,
\mathcal{L}_\alpha^\beta.
\end{equation}
Next, if as usual $\delta_{j_2}^{j_1}$ denotes the Kronecker symbol,
equal to $1$ if $j_1 = j_2$ and to $0$ otherwise, we can
shortly the various occuring partial derivatives of the monomial
$z^\alpha$ as:
\[
\aligned
\frac{\partial(z^\alpha)}{\partial z_{j_1}}
&
=
\alpha_{j_1}\,z^{\alpha-\epsilon_{j_1}},\,\,\,
\frac{\partial^2(z^\alpha)}{\partial z_{j_1}\partial z_{j_2}}
=
\alpha_{j_1}
\big(\alpha_{j_2}-\delta_{j_2}^{j_1}\big)\,
z^{\alpha-\epsilon_{j_1}-\epsilon_{j_2}},\,
\cdots\cdots,\,\,
\\
\frac{\partial^e(z^\alpha)}{\partial z_{j_1}\partial z_{j_2}\cdots
\partial z_{j_e}}
&
=
\alpha_{j_1}
\big(\alpha_{j_2}-\delta_{j_2}^{j_1}\big)
\cdots
\big(\alpha_{j_e}-\delta_{j_e}^{j_1}
-\cdots-
\delta_{j_e}^{j_{e-1}}\big)\,
z^{\alpha-\epsilon_{j_1}-\cdots-\epsilon_{j_e}}.
\endaligned
\] 
It follows that, for every $e \leqslant n$, the $(e + 1)$-th family of
equations, after equating to zero the coefficients of the monomial
$z^{\rho - \epsilon_{ j_1} - \cdots - \epsilon{ j_e}}$ and replacing
$\alpha$ by $\rho - \beta$, identifies to the collection:
\def\theequation{${\bf e}_{j_1j_2\cdots j_e\rho}$}\begin{equation}
\aligned
0
&
=
\sum_{\vert\beta\vert\leqslant n}\,
(\rho_{j_1}-\beta_{j_1})
\big(\rho_{j_2}-\beta_{j_2}-\delta_{j_2}^{j_1}\big)
\cdots
\big(\rho_{j_e}-\beta_{j_e}-\delta_{j_e}^{j_1}
-\cdots-
\delta_{j_e}^{j_{e-1}}\big)\,
\mathcal{L}_{\rho-\beta}^\beta
+
\\
&
\ \ \ \ \ \ \ \ \ \ \ \ \ \ \ \ \ \ \ \ \ \ \ \ \ \ \ \ \ \ 
+
{\sf R}_{j_1j_2\cdots j_e\rho}(a, \Lambda),
\endaligned
\end{equation}
where each second term ${\sf R}_{j_1j_2 \cdots j_e \rho} (a,
\Lambda)$, here considered as being just a remainder, is the
coefficient of $z^{\rho - \epsilon_{j_1} - \cdots - \epsilon_{ j_e}}$
in the second term of the equation $({\bf e}_{ j_1 j_2 \cdots j_e})$
and hence is clearly bilinear in $\big( a_\gamma, \Lambda_k^l\big)$.

Thus, we have written a constant coefficient system of linear
equations having the $\mathcal{ L}_{\rho - \beta}^\beta$ as unknowns,
$\vert \beta \vert \leqslant n$. As in~\cite{ pau2008, rou2007a}, we
now claim that the determinant of its matrix is nonzero.

Indeed, for each fixed multiindex $\rho$, the matrix whose column
$C_\beta$ consists of the partial derivatives of order at most $n$ of
the monomial $z^{ \rho - \beta}$ has the same determinant, at the
point $(1, 1, \dots, 1)$ as the linear subsystem $({\bf 0}_\rho)$, $(
{\bf 1}_{ j_1 \rho})$, \dots, $( {\bf e}_{ j_1 \cdots j_n \rho})$ we
want to solve, where $j_1, \dots, j_n = 1, \dots, n+1$.  Therefore, if
the determinant would be zero, we would by linear combination, derive
the existence of a {\em not} identically zero polynomial:
\[
Q(z)
:=
\sum_{\beta}\,c_\beta\,z^{\rho-\beta}
\]
all of whose partial derivatives of order $\leqslant n$ vanish at $(1,
\dots, 1)$. Hence the same would be true of:
\[
P(z)
:=
z^\rho\,
Q\big(
1/z_1, \dots, 1/z_{n+1}\big)
=
\sum_\beta\,c_\beta\,z^\beta,
\]
and this would imply $P \equiv 0$, in contradiction to the assumption.

Thus for each fixed $\rho$, Cramer's rule solves the system for the
$\mathcal{ L}_\alpha^\beta$ with $\alpha + \beta = \rho$, and the
solution is then obviously bilinear in $(a, \Lambda)$.
\endproof

\subsection*{ Invariance under reparametrization and logarithmic
versions} We would like to make two final remarks, useful in
applications. At first, similarly as it was pointed out in~\cite{
pau2008, rou2007a}, we claim that all vector fields constructed above
are invariant under the group ${\sf G}_n$ of $n$-jets at the origin of
local reparametrizations $\phi ( \zeta) = \zeta + \phi'' ( 0) \,
\frac{ \zeta^2}{ 2!} + \cdots + \phi^{ ( n) } ( 0) \, \frac{ \zeta^n
}{ n !} + \cdots$ of $(\C, 0) \ni \zeta$ that are tangent to the
identity\footnote{\, As a result, our two theorems can be applied in
the framework of Demailly-Semple jets (\cite{ dmr2008}).}, which acts
on the jets $(z_{ i_1}', z_{ i_2}'', z_{ i_3}''', \dots)$ by
transforming them to:
\[
w_{i_1}'
:=
z_{i_1}',
\ \ \ \ \ \ \ \ \ \ \ \
w_{i_2}''
:=
z_{i_2}''+\phi''\,z_{i_2}',
\ \ \ \ \ \ \ \ \ \ \ \
w_{i_3}'''
:=
z_{i_3}'''+3\phi''z_{i_3}''+\phi'''z_{i_3}',\,\,
\dots.
\]
Such a transformation makes a diffeomorphism of $J_{ \sf vert}^n (
\mathcal{ X})$, its inverse being associated to $\phi^{ -1} ( \zeta)$,
and we must verify that our 4 families of tangent vector fields ${\sf
T}_\alpha$, ${\sf T}_\alpha^{ \ell_1, \dots, \ell_{ n+1}}$, ${\sf
T}_i$ and ${\sf T}_\Lambda$ are left unchanged under this
diffeomorphism.  Indeed, the claim trivially holds true for both the
${\sf T}_\alpha^{ \ell_1, \dots, \ell_{ n+1}}$ and the ${\sf T}_i$,
because they incorporate absolutely no $z_{ i_1}', z_{ i_2}'', \dots,
z_{ i_n}^{ (n)}$. Next, the $\frac{ \partial}{ \partial a_\beta}$ in
the ${\sf T}_\alpha$ are clearly left unchanged, while their
coefficients are all, say in the case $n = 3$ to fix ideas, of the
Wronskian-like form:
\[
\footnotesize
\aligned
\left\vert
\begin{array}{ccc}
f' & g' & h' 
\\
f'' & g'' & h''
\\
f''' & g''' & h'''
\end{array}
\right\vert
\equiv
\left\vert
\begin{array}{ccc}
f' & g' & h' 
\\
f'' + \phi'' f' & g'' + \phi'' g' & h'' + \phi'' h'
\\
f''' + 3\phi'' f'' + \phi''' f' & 
g''' + 3\phi'' g'' + \phi''' g' &
h''' + 3\phi'' h'' + \phi''' h'
\end{array}
\right\vert,
\endaligned
\]
where $f, g, h \in \C[ z_1, z_2, z_3 ]$ are some polynomials, but then
such a determinant remains unchanged, thanks to obvious line
manipulations. The general case $n \geqslant 3$ is similar. Finally,
erasing indices and again for $n = 3$, we give the formal reason why
the ${\sf T}_\Lambda$ are also invariant. The transformation is $w' =
z'$, $w'' = z'' + \phi'' z'$, $w''' = z''' + 3 \phi'' z'' + \phi'''
z'$ and it replaces the basic vector fields by: 
\[
{\textstyle{\frac{\partial}{\partial z'}}}
=
{\textstyle{\frac{\partial}{\partial w'}}}
+
\phi''\,
{\textstyle{\frac{\partial}{\partial w''}}}
+
\phi'''\,
{\textstyle{\frac{\partial}{\partial w'''}}},
\ \ \ \ \ \
{\textstyle{\frac{\partial}{\partial z''}}}
=
{\textstyle{\frac{\partial}{\partial w''}}}
+
3\phi''\,
{\textstyle{\frac{\partial}{\partial w'''}}},
\ \ \ \ \ \
{\textstyle{\frac{\partial}{\partial z'''}}}
=
{\textstyle{\frac{\partial}{\partial w'''}}},
\]
so that $z' \frac{ \partial }{ \partial z'} + z'' \frac{ \partial
}{\partial z''} + z''' \frac{ \partial }{ \partial z'''} = w' \frac{
\partial }{ \partial w'} + w'' \frac{ \partial }{\partial w''} + w'''
\frac{ \partial }{ \partial w'''}$ is invariant.

Another argument (transmitted to 
us by Erwan Rousseau) for invariance under reparametrization would be to
say that the system of linear equations that the coefficients ${\sf
Z}_i$, ${\sf A}_\alpha$, ${\sf Z}_k'$, ${\sf Z}_k''$, \dots, ${\sf
Z}_k^{(n)}$ of a general tangent vector field ${\sf T}$ have to
satisfy:
\[
0
=
{\sf T}\big[
{\textstyle{\sum_\alpha}}\,a_\alpha z^\alpha\big]
=
{\sf T}\big[
{\textstyle{\sum_\alpha}}\,a_\alpha (z^\alpha)'\big]
=
{\sf T}\big[
{\textstyle{\sum_\alpha}}\,a_\alpha (z^\alpha)''\big]
=
{\sf T}\big[
{\textstyle{\sum_\alpha}}\,a_\alpha (z^\alpha)'''\big]
=\cdots,
\]
is transformed, after reparametrization, into a system:
\[
\aligned
0
&
=
{\sf T}\big[
{\textstyle{\sum_\alpha}}\,a_\alpha z^\alpha\big]
=
{\sf T}\big[
{\textstyle{\sum_\alpha}}\,a_\alpha (z^\alpha)'\big]
=
{\sf T}\big[
{\textstyle{\sum_\alpha}}\,a_\alpha (z^\alpha)''
+
\phi''
{\textstyle{\sum_\alpha}}\,a_\alpha (z^\alpha)'
\big]
\\
0
&
=
{\sf T}\big[
{\textstyle{\sum_\alpha}}\,a_\alpha (z^\alpha)'''
+
3\phi''\,
{\textstyle{\sum_\alpha}}\,a_\alpha (z^\alpha)''
+
\phi'''
{\textstyle{\sum_\alpha}}\,a_\alpha (z^\alpha)'
\big]
=\cdots
\endaligned
\]
which is completely equivalent to the first one, thanks to obvious
linear combinations, so that any solution to this linear system is
{\em a priori} forced to be invariant. 

\smallskip

The second remark is that one may adapt the formalism provided here to
show that the global generation property holds in a logarithmic
setting with the {\em same} specific pole orders $c_n = \frac{ n^2 +
5n}{ 2}$ ({\em cf.}~\cite{ rou2007b} for $n = 3$) or $c_n = n^2 +
2n$. Application to effective algebraic degeneracy of entire
holomorphic maps in the complement of a generic hypersurface $X
\subset \P^{ n+1} ( \C)$ are therefore also possible.

\section*{ \S5.~Appendix: a determinantal identity}
\label{Section-5}

\subsection*{ Proof of the combinatorial lemma}
At the beginning of Section~3, a determinant left aside
had to be computed. We drop the index $i$ and we denote it shortly by:
\[
\small
\aligned
\Delta
:=
\left\vert
\begin{array}{cccc}
z' & (z^2)' & \cdots & (z^n)'
\\
z'' & (z^2)'' & \cdots & (z^n)''
\\
\cdot\cdot & \cdot\cdot & \cdots & \cdot\cdot
\\
z^{(n)} & (z^2)^{(n)} & \cdots & (z^n)^{(n)}
\end{array}
\right\vert.
\endaligned
\]
On the first line, the entry of the $k$-th column is $(z^k)' = k
z^{k-1} z'$. The trick is then to write the entry of the second line
inside the same column as $k (z^{k-1} z')'$, {\em etc.}, and generally
the entry of the $\kappa$-th line as $k (z^{ k-1} z')^{ (\kappa-1)}$, 
so that:
\[
\footnotesize
\aligned
\Delta
=
\left\vert
\begin{array}{ccccc}
z' & 2zz' & 3z^2z' & \cdots & nz^{n-1}z'
\\
z'' & 2(zz')' & 3(z^2z')' & \cdots & n(z^{n-1}z')'
\\
z''' & 2(zz')'' & 3(z^2z')'' & \cdots & n(z^{n-1}z')''
\\
\cdot\cdot & \cdot\cdot & \cdot\cdot & \cdots & \cdot\cdot
\\
z^{(n)} & 2(zz')^{(n-1)} & 3(z^2z')^{(n-1)} & \cdots &
n(z^{n-1}z')^{(n-1)}
\end{array}
\right\vert.
\endaligned
\]
We see that $2 \cdot 3 \cdots n$ from the columns comes into factor.
Next, using Leibniz's formula for the derivative of a product, we may
expand $(z^{ k-1} z')^{ (\kappa - 1)}$ just as $\sum_{ 0 \leqslant
\lambda_1 \leqslant \kappa - 1} \binom{ \kappa-1}{ \lambda_1} (z^{
k-1})^{(\lambda_1)} z^{(1+\kappa-1-\lambda_1)}$.  By subtracting to
the $k$-th column the first one multiplied by $z^{ k-1}$, the
$(\kappa, k)$-entry then becomes $\sum_{ 1 \leqslant \lambda_1
\leqslant \kappa - 1} \binom{ \kappa-1}{ \lambda_1} (z^{
k-1})^{(\lambda_1)} z^{(\kappa-\lambda_1)}$, where now the sum starts
from $\lambda_1 = 1$. In particular, the $(1,2)$-, $(1, 3)$-, \dots,
$(0, n)$-entries all become null. By expanding the determinant along
its first line, we are therefore left with an $(n-1)\times (n-1)$
determinant (notice the necessary shift of indices):
\[
\small
\aligned
\frac{\Delta}{n!z'}
=
\bigg\vert
\sum_{1\leqslant\lambda_1\leqslant\kappa}
{\textstyle{\binom{\kappa}{\lambda_1}}}(z^k)^{(\lambda_1)}
z^{(1+\kappa-\lambda_1)}
\bigg\vert_{1\leqslant\kappa\leqslant n-1}^{
1\leqslant k\leqslant n-1}.
\endaligned
\]
Iterating the trick, we again write:
\[
\small
\aligned
(z^k)^{(\lambda_1)}
=
k(z^{k-1}z')^{(\lambda_1-1)}
=
k\,
{\textstyle{\sum_{0\leqslant\lambda_2\leqslant\lambda_1-1}}}\,
{\textstyle{\binom{\lambda_1-1}{\lambda_2}}}\,
(z^{k-1})^{(\lambda_2)}\,
z^{(1+\lambda_1-1-\lambda_2)}.
\endaligned
\]
We again see that $2 \cdot 3 \cdots (n-1)$ from the columns comes into
factor, and then substituting the computed value of $(z^k)^{
(\lambda_1)}$, we get:
\[
\footnotesize
\aligned
\frac{\Delta}{n!z'}
=
(n-1)!\cdot
\bigg\vert
\sum_{1\leqslant\lambda_1\leqslant\kappa}\,
&
\sum_{0\leqslant\lambda_2\leqslant\lambda_1-1}\,
{\textstyle{\binom{\kappa}{\lambda_1}}}\,
{\textstyle{\binom{\lambda_1-1}{\lambda_2}}}\cdot
\\
&
\cdot(z^{k-1})^{(\lambda_2)}\,
z^{(1+\kappa-\lambda_1)}\,
z^{(\lambda_1-\lambda_2)}
\bigg\vert_{1\leqslant\kappa\leqslant n-1}^{
1\leqslant k\leqslant n-1}.
\endaligned
\]
The $(\kappa, 1)$-entry inside the first colum is equal to $\sum_{ 1
\leqslant \lambda_1 \leqslant \kappa}\, \binom{ \kappa}{ \lambda_1} \,
z^{ ( 1 + \kappa - \lambda_1)}\, z^{ (\lambda_1)}$, because the terms
$(z^0)^{ ( \lambda_2)}$ with $\lambda_2 \geqslant 1$ are null. By
subtracting to the $k$-th column the first one multiplied by $z^{
k-1}$, the $(\kappa, k)$-th entry written above is slightly modified:
the sum involving $\lambda_2$ is then just replaced by $\sum_{ 1
\leqslant \lambda_2 \leqslant \lambda_1 - 1}$. Moreover, the $(1,
2)$-, $(1, 3)$-, \dots, $(1, n-1)$- entries all become null, while the
$(1, 1)$ entry is $\binom{ 1}{ 1} z' z'$. By expanding the determinant
along its first line, we are therefore left with an $(n-2) \times (
n-2)$ determinant:
\[
\footnotesize
\aligned
\frac{\Delta}{n!(n-1)!\,z'(z')^2}
=
\bigg\vert
\sum_{1\leqslant\lambda_1\leqslant\kappa}\,
&
\sum_{1\leqslant\lambda_2\leqslant\lambda_1-1}\,
{\textstyle{\binom{\kappa}{\lambda_1}}}\,
{\textstyle{\binom{\lambda_1-1}{\lambda_2}}}\cdot
\\
&
\cdot
(z^{k-1})^{(\lambda_2)}\,
z^{(1+\kappa-\lambda_1)}\,
z^{(\lambda_1-\lambda_2)}
\bigg\vert_{2\leqslant\kappa\leqslant n-1}^{2\leqslant k\leqslant n-1}.
\endaligned
\]
We now have to change the indices. We at first set $k' := k - 1$ and
$\kappa ' := \kappa - 1$ and the determinant just obtained becomes:
\[
\footnotesize
\aligned
\bigg\vert
\sum_{1\leqslant\lambda_1\leqslant\kappa'+1}\,
&
\sum_{1\leqslant\lambda_2\leqslant\lambda_1-1}\,
{\textstyle{\binom{\kappa'+1}{\lambda_1}}}\,
{\textstyle{\binom{\lambda_1-1}{\lambda_2}}}\cdot
\\
&
\cdot
(z^{k'})^{(\lambda_2)}\,
z^{(2+\kappa'-\lambda_1)}\,
z^{(\lambda_1-\lambda_2)}
\bigg\vert_{1\leqslant\kappa'\leqslant n-2}^{
1\leqslant k'\leqslant n-2}.
\endaligned
\]
Next, if we set $\lambda_1 ' := \lambda_1 - 1$ and if we observe the
identification of sums:
\[
\footnotesize
\aligned
\sum_{1\leqslant\lambda_1\leqslant\kappa'+1}\,
\sum_{1\leqslant\lambda_2\leqslant\lambda_1-1}
(\bullet)
=
\sum_{0\leqslant\lambda_1'\leqslant\kappa'}\,
\sum_{1\leqslant\lambda_2\leqslant\lambda_1'}
(\bullet)
=
\sum_{1\leqslant\lambda_2\leqslant\lambda_1'\leqslant\kappa'}
(\bullet),
\endaligned
\]
then our determinant simply becomes, after erasing the primes: 
\[
\footnotesize
\aligned
\bigg\vert
\sum_{1\leqslant\lambda_2\leqslant\lambda_1\leqslant\kappa}
{\textstyle{\binom{\kappa+1}{\lambda_1+1}}}\,
{\textstyle{\binom{\lambda_1}{\lambda_2}}}\cdot
z^{(\kappa-\lambda_1+1)}\,
z^{(\lambda_1-\lambda_2+1)}\cdot
(z^{k})^{(\lambda_2)}
\bigg\vert_{1\leqslant\kappa\leqslant n-2}^{
1\leqslant k\leqslant n-2}.
\endaligned
\]
Performing the same computational and transformational processes, the
result of the next step will be:
\[
\footnotesize
\aligned
\frac{\Delta}{n!(n-1)!(n-2)!\,z'(z')^2(z')^3}
&
=
\bigg\vert
\sum_{1\leqslant\lambda_3\leqslant\lambda_2\leqslant\lambda_1
\leqslant\kappa}
{\textstyle{\binom{\kappa+2}{\lambda_1+2}}}\,
{\textstyle{\binom{\lambda_1+1}{\lambda_2+1}}}\,
{\textstyle{\binom{\lambda_2}{\lambda_3}}}\cdot
\\
&
\ \ \ \ \
\cdot
z^{(\kappa-\lambda_1+1)}\,
z^{(\lambda_1-\lambda_2+1)}\,
z^{(\lambda_2-\lambda_3+1)}\cdot
(z^{k})^{(\lambda_3)}
\bigg\vert_{1\leqslant\kappa\leqslant n-3}^{
1\leqslant k\leqslant n-3}.
\endaligned
\]
The induction is now clear, and at the end one obtains a $1 \times 1$
determinant $\big\vert (\bullet) \big\vert_{ 1\leqslant \kappa
\leqslant 1}^{ 1 \leqslant k \leqslant 1}$ with a sum $\sum_{ 1
\leqslant \lambda_{ n -1} \leqslant \cdots \leqslant \lambda_1
\leqslant \kappa} ( \bullet )$ inside which necessarily $k = \kappa
= \lambda_1 = \cdots = \lambda_{ n-1}$, so that this last $1 \times 1$
determinant equals:
\[
{\textstyle{\binom{n-1}{n-1}}}\cdots
{\textstyle{\binom{1}{1}}}\,
z'\cdots z' \cdot z'
=
1!\,(z')^n,
\]
and this final observation completes the
proof of the combinatorial lemma.
\qed

\vfill\end{document}